\documentclass[final,onefignum,onetabnum]{siamart190516}

\usepackage{algorithm}
\usepackage{algpseudocode}
\usepackage{graphicx, epsfig}
\usepackage{latexsym}
\usepackage{epsfig}
\usepackage{latexsym}
\usepackage{amstext}
\usepackage{ulem}
\usepackage{caption}
\usepackage{amssymb}
\usepackage{array}
\usepackage{dcolumn}
\usepackage{amsmath, mathtools}
\usepackage{subcaption}
\usepackage{epsfig, graphicx}
\usepackage{color}
\usepackage{amsfonts,bm,color}
\usepackage{verbatim}
\usepackage{float}
\usepackage{latexsym}
\usepackage[svgnames]{xcolor}
\usepackage{colortbl}
\usepackage{float}
\usepackage{wrapfig}
\usepackage{multirow}
 \usepackage{stackrel}
 \usepackage{changepage}
\newcommand{\beq}{\begin{eqnarray}}
\newcommand{\eeq}{\end{eqnarray}}
\newcommand{\beqq}{\begin{eqnarray*}}
\newcommand{\eeqq}{\end{eqnarray*}}

\newcommand{\eps}{\varepsilon}

\newcommand{\ds}{\displaystyle}

\title{Accelerated Simulation Algorithms for Extreme First-Passage Problems with General Emission Profiles}
\author{E. Akame Mfoumou\thanks{\'Ecole Normale Sup\'erieure, Universit\'e PSL, Paris, France.}
\and D. Holcman\thanks{\'Ecole Normale Sup\'erieure, Universit\'e PSL, Paris, France and Churchill College, CB3 0DS, United Kingdom.}}
\date{\today}
\begin{document}

\maketitle
\begin{abstract}
Fastest arrival events, where the first among many diffusing particles reaches a target, are central in triggering signal initiation in molecular stochastic systems. Classical approaches to simulate such events rely on full trajectory generation of all particles, leading to prohibitive computational costs in the large particle number regime. In this work, we present a general simulation framework for efficiently generating order statistics of arrival times by exploiting asymptotic first-passage distributions.  { This framework applies to diffusion processes in bounded domains with localized absorbing targets, for which short-time first-passage asymptotics are available, such as Brownian motion in dimensions one, two, and three.} Starting with the case of instantaneous emission, we derive and implement a recursive inverse transform algorithm to simulate the first $k$ arrivals without tracking particle trajectories. We  extend this algorithm to time-dependent emission profiles via { an iterative approach, enabling the simulation of extreme statistics in systems with temporal injection, ranging from rapid to prolonged emission}. Additionally, we provide asymptotic estimates of the mean fastest arrival time. To conclude, the present acceleration algorithm which bypasses Brownian simulations of trajectories can be used for spatial reaction networks, rare event detection, or diffusion-controlled activation.
\end{abstract}
\section{Introduction}
Extreme first-passage events, in which the fastest particle among a large population reaches a target, play a key role in many stochastic processes such as molecular activation, synaptic signaling, gene regulation, telomere length or fertility \cite{reynaud2015so,Sokolov2019extreme,schuss2018cells,schehr2014exact,Basnayake2018}. In these examples, the dynamic is dominated not by the average arrival time but by the leading edge of the distribution — the earliest arrivals that can trigger downstream events. The study of such rare events has led to the development of extreme value theory for Brownian motion \cite{kotz2000extreme,coles2001introduction,leadbetter2012extremes}, often coupled with the narrow escape framework to account for geometric constraints \cite{Schuss:Book,holcman2015stochastic}.  However, simulating these fastest arrival times remains a challenge for large numbers. Indeed,  standard Brownian dynamic requires tracking all particle trajectories \cite{Basnayake2018}, which is computationally inefficient, especially in high-copy-number or multi-target settings. Recent advances have introduced asymptotic formulas for the fastest arrival time distributions and their associated splitting probabilities, enabling analytical progress \cite{Basnayake2018,Holcmanschuss2018,basnayake2020extremecusp,lawley2024hittingfast,lawley2020distribution}. { In particular, explicit formulas \cite{lawley2020distribution} use joint distribution and moments of the first $k$ arrival times of independent diffusing particles, including asymptotic expressions involving the Lambert $W$ function. While these results provide a general probabilistic description of extreme arrival statistics, they are primarily analytical and do not directly yield efficient simulation procedures. Thus practical simulation tools that leverage these formulas to replace full-path simulations are still under development.
We present here a simulation approach that relies on explicit numerical algorithms: rather than sampling from closed-form distributions, we will develop a recursive algorithm based on asymptotic survival laws that enables direct simulation of the first $k$ arrivals without generating particle trajectories. { This approach extends  to multi-target geometries, time-dependent emission profiles, and killing processes.   We shall use here explicit asymptotic formulas and the construction of efficient recursive simulation algorithms for extreme first-passage statistics, parallel to previous approaches \cite{grebenkov2025,lawley2020distribution}.\\
The present recursive construction differs from a direct inverse CDF sampling of known extreme-value distributions, as it relies only on local survival probabilities and does not require explicit knowledge or inversion of the full joint distribution. Finally, we will see here that the principal advantage of the framework lies in its scalability: the runtime scales linearly with the number of arrivals $k$ to be sampled, but is essentially independent of the total number of particles $N$. } \\
We recall that numerical approaches to reduce computational cost have emerged to bridge large particle baths with small target domains \cite{nadler2001stationary}, particularly in settings where discrete stochastic events are embedded within a continuum framework \cite{erban2007reactive,singer2008partially,erban2009stochastic}. These hybrid schemes have been successfully applied in cell signaling models to convert spatially averaged concentrations, such as mean calcium levels, into discrete particle representations for simulating localized events such as triggering synaptic transmission from discrete calcium binding events \cite{guerrier2016hybrid,guerrier2018first}. Stochastic frameworks that account for partial absorption, reactive boundaries, and finite reaction windows in diffusion-limited systems \cite{erban2007reactive,erban2009stochastic} are particularly relevant when particles face competition between target binding and decay, leading to corrections in classical narrow escape statistics \cite{Holcmanschuss2018}.\\
Simulating stochastic and anomalous transport processes in biology using multiscale numerical integration \cite{burrage2004foundations} has improved many current simulation strategies, especially when exact trajectories are expensive to compute. This is the case for fractional-order models in biological systems \cite{burrage2024fractional} where accurate solvers have been designed for fractional differential equations \cite{brugnano2024spectrally}, and efficient numerical schemes for diffusion on curved geometries, such as spherical surfaces in 3D \cite{burrage2022effective}.  These developments highlight the importance of combining analytical insight with algorithmic accuracy to simulate rare events. We will show how discrete outcomes (binding, death, or internalization) can be statistically captured through specific  designed sampling methods. We will extend the notion of accelerated simulation for extreme first-passage events to incorporate the time distribution of injection and exponentially distributed killing, capturing the interplay between mortality and spatial escape in a unified and computationally efficient framework. For that goal, we introduce a simulation strategy for extreme first-passage problems that directly samples arrival times and targets without computing trajectories.\\
{ The present framework is designed for stochastic processes whose single-particle first-passage time distribution is either known explicitly or admits a short-time asymptotic approximations. Our primary focus is on Brownian diffusion in bounded domains with small absorbing regions, where such asymptotics have been established. The algorithm may also be applicable to more general diffusion processes—including those in complex potentials or with spatially varying diffusivity—provided that analogous survival probability expansions can be obtained. More broadly, the present approach could be extended to stochastic and numerical approaches for transport and flows in complex media, including porous and electrokinetic systems \cite{RefIJFT2025,Ahmad2025IJFT}.}}\\
We begin in Section 1 with the case of instantaneous emission, where all particles are released simultaneously. Using asymptotics for first-passage time distribution for the fastest, we derive recursive inverse transform algorithms for simulating the first $k-$arrivals, incorporating splitting probabilities for multiple targets. In Section 2, we extend the framework to general emission profiles — including bell-shape or exponentially modulated injections — by using successive convolutions for effective survival probabilities. We show that the simulation of fastest arrivals under such general emission laws can still be accomplished via an accelerated scheme.\\
Throughout the paper, we validate our algorithms against classical Brownian simulations and provide asymptotic expansions for key quantities such as Mean First Arrival Time (MFAT). These expansions clarify the distinct scaling regimes induced by the emission rate relative to the diffusive timescale. In particular, we recover and extend results involving the Lambert $W-$function in the characterization of MFATs. Our results offer a general and efficient framework for simulating fastest particle driven events in two and three dimensional domains with localized targets and under time-dependent emission dynamics. { Applications range from cell biology and materials science, to chemical physics, wherever rare but rapid activation is essential. Indeed, time-dependent injection processes naturally arise in a wide range of biological systems where activation is triggered by rare but fast events. For example, neurotransmitter release at synapses is governed by stochastic vesicle fusion events controlled by calcium binding, while calcium signaling in dendritic spines involves transient fluxes that activate Ryanodine receptors (RyRs) through the arrival of the fastest ions \cite{basnayake2019fastest}. Similarly, time-dependent ATP production by mitochondria leads to delayed diffusion toward molecular targets, while sensory transduction processes, such as in Drosophila photoreceptors, involves molecular cascades followed by diffusion in confined geometries, where activation is controlled by the fastest signaling molecules. These examples illustrate that extreme first-passage events under temporally structured injection are central to molecular activation mechanisms in cells \cite{fain2019sensory}}.
\section{Instantaneous Emission: Recursive Simulation Algorithm}
\subsection{Motivation and Description}
In the setting of instantaneous emission, all $n$ particles are released simultaneously from a point source $y$ within a bounded domain. The goal is to simulate the arrival time statistics of the fastest (and also successive) particles reaching one or more absorbing targets on the boundary. Rather than simulating all particle trajectories and tracking the first arrival, we instead here used known asymptotic formulas for the distribution of first arrival times $F(t)$ and their corresponding order statistics. { This approach is relevant when short-time asymptotic expressions for the survival probability of a single particle are known.  This is for example the case for Brownian diffusion in a bounded domain with small absorbing regions, as in the classical narrow escape framework \cite{holcman2015stochastic}. In this regime, the geometry of the domain enters through a finite set of geodesic distances to the targets that can be directly exploited in the sampling algorithm \cite{basnayake2020extremecusp}.}\\
We start by recalling that the survival probability $S(t)$  of a single Brownian particle serves to express $S(t)^n$ of all $n$ particles that survive up to time $t$, and thus the distribution function of the fastest arrival time is $F_n(t) = 1 - S(t)^n$. Indeed, \cite{Holcmanschuss2018}
\beq
 Pr\lbrace\tau^n >t\rbrace = Pr\lbrace\tau >t\rbrace^n = S(t)^n.
\eeq
We will use below the inverse of this distribution to sample the arrival time of the fastest particle. When there is only one small absorbing region for a Brownian motion characterized by a diffusion coefficient $D$, the asymptotic expansion for survival probability $S(t)$ is given by \cite{Basnayake2018}
\beq\label{asympexpressions}
    S(t) \approx
    \begin{cases}
        1-\frac{\sqrt{4Dt}}{\delta \sqrt{\pi}} \exp\left(- \frac{\delta^2}{4Dt}\right) \qquad \text{(1-Dim)} \\
        1-\frac{\sqrt{2}\pi Dt}{2\log\left(\frac{1}{\eps}\right)\delta^2}\exp\left(-\frac{\delta^2}{4Dt}\right) \qquad \text{(2-Dim)} \\
        1-\frac{a^2}{\delta\sqrt{\pi Dt}}\exp\left(-\frac{\delta^2}{4Dt}\right) \,\,\,  \qquad \text{(3-Dim)}
    \end{cases}
\eeq
where $\delta$ is the geodesic distance between the source and the center of the absorbing patch of size $2\eps$ (2D) and $a$ the radius of the exit in three dimension.  When there are $M$ small, well-separated absorbing sites located at distances  $\delta_1,\ \delta_2, ... , \delta_m$ from the source in a domain, the total survival probability is the product of the individual survival probabilities near each target. {In one dimension, we have M=2. The case of multiple target ($M>2$) absorbing sites could however be relevant on a graph with multiple ends \cite{dora2020active}. The first arrival is dominated by the nearest ones. In general, the total survival probability which is product of individual one can be approximated as follows (see Appendix A):}
\beq
    S(t) \approx
    \begin{cases}
        1 - \sum\limits_i^M \frac{\sqrt{4Dt}}{\delta_i \sqrt{\pi}} \exp\left(- \frac{\delta_i^2}{4Dt}\right) \\
        1-\sum\limits_i^M\frac{\sqrt{2}\pi Dt}{2\log\left(\frac{1}{\eps_i}\right)\delta_i^2}\exp\left(-\frac{\delta^2_i}{4Dt}\right) \label{survival123}\\
        1-\sum\limits_i^M\frac{a_i^2}{\delta_i\sqrt{\pi Dt}}\exp\left(-\frac{\delta^2_i}{4Dt}\right)
    \end{cases}
\eeq
The survival probability of the fastest is given by the expression
\beq\label{distrbutiontaun}
    Pr\lbrace\tau^n >t\rbrace \approx
    \begin{cases}
        \exp\left(-n\sum\limits_i^M \frac{\sqrt{4Dt}}{\delta_i \sqrt{\pi}} \exp\left(- \frac{\delta_i^2}{4Dt}\right) \right) \\
        \exp\left(-n\sum\limits_i^M\frac{\sqrt{2}\pi Dt}{2\log\left(\frac{1}{\eps_i}\right)\delta_i^2}\exp\left(-\frac{\delta^2_i}{4Dt}\right)\right) \\
        \exp\left(-n\sum\limits_i^M\frac{a_i^2}{\delta_i\sqrt{\pi Dt}}\exp\left(-\frac{\delta^2_i}{4Dt}\right)\right).
    \end{cases}
\eeq
We now describe the heuristic for an accelerated simulation algorithm: to sample the first $k$ arrivals, we use a recursive construction based on the known order statistics for independent identically distributed (i.i.d.) random variables. The cumulative distribution $F(t) = 1 - S(t)$ function (CDF) of the first-passage time for a single particle allows computing the conditional distribution of the $(k+1)$-th arrival time given that $k$ arrivals have occurred by time $t_k$ (see derivation in the Appendix)
\beq
\mathbb{P}(\tau_{k+1} > t \mid \tau_k = t_k) = \left(\frac{S(t)}{S(t_k)}\right)^{n-k}\approx \exp (-(n-k)(F(t)-F(t_k))).
\eeq
In practice, given a uniform random variable $\mathbb{P}(\tau_{k+1} > t \mid \tau_k = t_k)=U_{k+1} \sim \mathcal{U}[0,1]$, the next arrival time $t_{k+1}$ satisfies:
\beq
F(t_{k+1}) = F(t_k) + \frac{\log(1/U_{k+1})}{n-k}.
\eeq
This leads to the following recursive algorithm:
\begin{itemize}
\item Set $t_0 = 0$ and initialize $F(t_0) = 0$.
\item For $k = 1, \dots, K$:
\begin{enumerate}
\item Draw $U_k \sim \mathcal{U}[0,1]$.
\item Compute $F(t_k) = F(t_{k-1}) + \frac{\log(1/U_k)}{n - (k - 1)}$.
\item Invert $F$ to find $t_k = F^{-1}(F(t_k))$.
\end{enumerate}
\end{itemize}
This method will allow us to simulate the first $k-$arrival times without ever simulating Brownian trajectories when the geodesic distances are known.
\subsection{Summary of Convergence and Validity}
The convergence of the above algorithm follows from the properties of order statistics. The distribution of the $k$-th order statistic among $n$ i.i.d. variables with cumulative distribution function $F$ is given by \cite{weiss1983order}:
\beq
f_{\tau_k^n}(t) = k \binom{n}{k} F(t)^{k-1} (1 - F(t))^{n - k} f(t),
\eeq
\label{eq:order}
where $f(t)$ is the density corresponding to $F(t)$. We define $F(t)$ via the asymptotic first-passage time distribution and using the exponential representation of order statistics. Moreover, the expression for the conditional distribution
\beq
\mathbb{P}(\tau_{k+1} > t \mid \tau_k = t_k) \approx \exp\left( -(n-k)(F(t) - F(t_k)) \right),
\eeq
shows that the increments of the CDF between arrivals behave like independent exponential variables. This algorithm avoids possible numerical instabilities of directly sampling from steep tail distributions and provides rapid convergence when $n$ is large and $F(t)$ captures the dominant short-time arrival behavior. In the next section, we provide the details of the algorithm in the three dimensions.
\subsection{Geometric Target Selection and Sampling Algorithm}
\begin{figure}[h!]
\begin{adjustwidth}{-2cm}{-2cm}
\centering
\includegraphics[scale = 0.75]{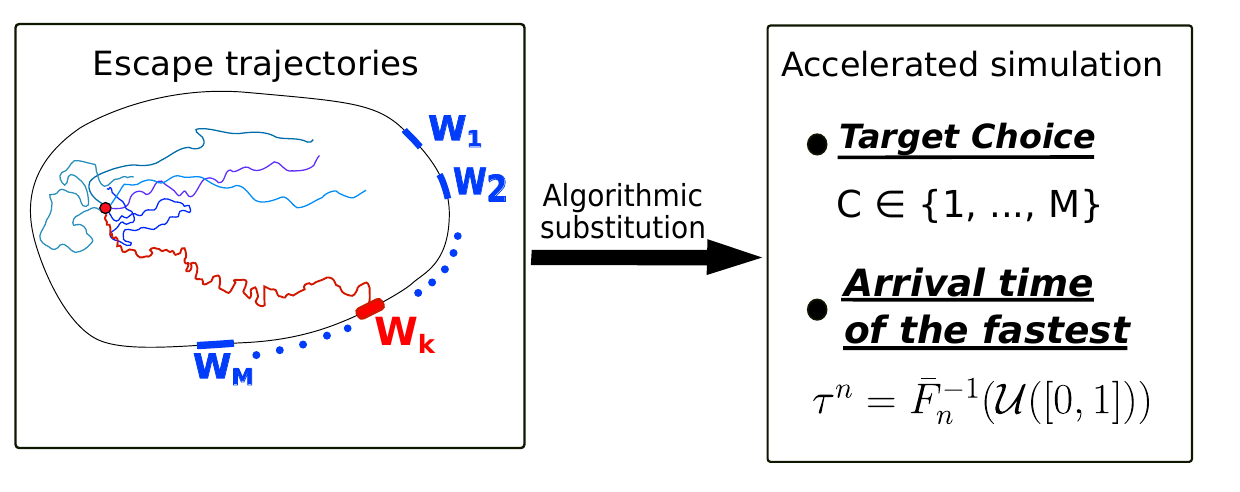}
\caption{\textbf{Substituting Brownian Simulations with Extreme Escape Statistics}. This schematic illustrates the core idea of the algorithm: replacing the computationally intensive simulation of Brownian trajectories up to the first escape event with the independent sampling of two random variables — the arrival time $\tau$ and the escape target $C$ — from precomputed probability distributions.}
\label{scheme}
\end{adjustwidth}
\end{figure}

The recursive sampling algorithm for fastest arrival times, as described above, will be able to handle the temporal statistics of successive escape. However, in systems with multiple absorbing targets, determining which target is reached first is also important. This step can be computed using the splitting probability, that is, the probability that a Brownian particle reaches a particular target among several \cite{karlin1981,Schuss:Book,schuss2019redundancy}. In this section, we account for spatial distribution of targets. Rather than simulating Brownian trajectories through space, we sample directly from known asymptotic expressions for escape probabilities and splitting ratios. This approach replaces Brownian simulations to the efficient sampling of two random variables: the escape time $\tau$ and the corresponding target index $C$ (fig. \ref{scheme}).
\paragraph{\bf Algorithm Description}
To simulate the earliest arrival among $n$ independent Brownian particles released from a common source $y$ in a bounded domain $\Omega$ containing $M$ absorbing targets, we propose the following algorithm:
\begin{enumerate}
\item \textbf{Step 1: Preprocess Geometry.} Compute the geodesic distances $\delta_i$ from source $y$ to each target $i = 1, \dots, M$ and sort them.
\item \textbf{Step 2: Sample Arrival Time.} Draw $U_1 \sim \mathcal{U}[0,1]$ and set:
\begin{equation}
    \tau = \bar{F}_n^{-1}(U_1),
\end{equation}
where $\bar{F}_n(t) = \mathbb{P}(\tau^n > t)$ is the survival function for the fastest particle among $n$.
\item \textbf{Step 3: Sample Target Identity.} Compute splitting probabilities $p_i$ using asymptotic expansions depending on the dimension (see below). Let $P_i = \sum_{j=1}^i p_j$, and draw $U_2 \sim \mathcal{U}[0,1]$. Find $i_0$ such that:
\begin{equation}
    P_{i_0-1} < U_2 < P_{i_0},
\end{equation}
and set the selected target index to $i_0$.
\end{enumerate}
We will recall below splitting probability expressions that depend on distance ratios $\lambda_{i,j}$ and geometric correction terms $k_{i,j}$ defined
\begin{equation}
\lambda_{i,j} = \begin{cases}
\ds \frac{\delta_i}{\delta_j} & \text{in 1D or 3D}, \\ \\
\ds \frac{\delta_i^2}{\delta_j^2} & \text{in 2D},
\end{cases} \quad
k_{i,j} = \begin{cases}
\ds \frac{\log(1/\varepsilon_i)}{\log(1/\varepsilon_j)} & \text{in 2D}, \\ \\
\ds \frac{a_j^2}{a_i^2} & \text{in 3D},
\end{cases}
\end{equation}
where $\varepsilon_i$ is the window width (2D) or $a_i$ is the window radius (3D). A pseudocode version of this algorithm is:
\begin{algorithm}[http!]
\caption{Extreme Statistics Simulation with Multiple Targets}
\footnotesize
\begin{algorithmic}[1]
\Require targets: list of $M$ targets
\Require $y$: source point
\Require $n$: number of particles
\State Compute and sort geodesic distances $\delta_1, \dots, \delta_M$
\State Draw $U_1 \sim \mathcal{U}[0,1]$; compute $\tau = \bar{F}_n^{-1}(U_1)$
\State Compute splitting probabilities $p_1, \dots, p_M$ using geometric asymptotics
\State Draw $U_2 \sim \mathcal{U}[0,1]$; assign target $i_0$ via cumulative $P_i$
\State \Return $(\tau, i_0)$
\end{algorithmic}
\end{algorithm}
\paragraph{Explicit Inversion Using Lambert $W$ Function}
When the dominant contribution to $\bar{F}_n(t)$ arises from the single closest target, we will derive below a dimension-specific inverses for the first arrival time $\tau = \bar{F}_n^{-1}(U)$. Indeed, when the dominant contribution to $\bar{F}_n(t)$ comes from the smallest distance $\delta_i$, an explicit formula for its inverse uses the Lambert $W$ function, depending on the spatial dimension:
\begin{equation}
    \bar{F}_n^{-1}(t) =
    \begin{cases}
    \ds \frac{\delta^2}{2DW_0\left( \ds \frac{2n^2}{\pi\ln(t)^2} \right)} \hbox{  dim } 1\\\\
    \ds \frac{\delta^2}{4DW_0\left( \ds \frac{\sqrt{2}\pi n}{8\ln\left(\frac{1}{\eps}\right) \ln\left(\frac{1}{t}\right)} \right)} \hbox{  dim } 2\\ \\
   \ds  -\frac{\delta^2}{2DW_{-1}\left(- \frac{\ds\pi\delta^4\ln \left(\frac{1}{t}\right)^2}{\displaystyle 2n^2a^4} \right)} \hbox{  dim } 3,
    \end{cases}
\end{equation}
where $W_0$ and $W{-1}$ are the principal and lower branches of the Lambert $W$ function, respectively,  defined by  $W(x)e^{W(x)} = x$.
\section{Explicit Algorithm Construction and Simulation of the $k$-th Fastest Arrival}
Building on the recursive algorithm developed for simulating the first arrival event, we now extend the method to generate the full sequence of the first $k$ arrival times among $n$ Brownian particles. This extension provides a practical and efficient alternative to simulating entire Brownian trajectories and enables precise modeling of sequential activation events in stochastic systems with many particles and small targets.
\subsection{Distribution of the First $k$ Arrival Times}
We derive here the statistical distribution of the ordered arrival times $\tau^n_1 < \tau^n_2 < \cdots < \tau^n_k$ of $n$ Brownian particles, that we will use to design our accelerated algorithm. We will compute the transition between the $k$-th and $(k+1)$-th arrivals.\\
Starting with the survival probability $S(t)$  of a single Brownian particle up to time $t$, the conditional probability that the $(k+1)$-th particle arrives later than time $t$, given that the $k$-th particle arrived at time $s$ is
\begin{equation}
    \mathbb{P}(\tau^n_{k+1} > t \mid \tau^n_k = s) = \left( \frac{S(t)}{S(s)} \right)^{n-k}, \quad \text{for } t > s.
\end{equation}
Using the short-time asymptotics for the survival probability $S(t)$ of a Brownian motion (see relation \ref{survival123})  $S(t) = 1 - F(t)$, where $F$ is the exit probability of the escape time of one Brownian,  we have $S(t)^n \approx \exp(-nF(t))$ so that the conditional survival probability becomes:
\beq
    Pr\lbrace\tau_{k+1}^n \geq t | \tau_k = s\rbrace \approx \exp[-(n-k)(F(t)-F(s))].
\eeq
Taking the logarithm and differentiating yields the recursive update rule for the exit probability  using the inverse transform method. Indeed, with $t_0=0$ and $U_{k+1} \sim \mathcal{U}(0,1)$, the arrival time of the $(k+1)$-th particle is computed by:
\beq\label{recursion}
    F(t_{k+1}) &= F(t_k)+\frac{\log\left(\frac{1}{U_{k+1}}\right)}{n-k}.
\eeq
We can also compute the joint distribution of the $k$-th and $(k+1)$-th arrival times $(\tau_k^n, \tau_{k+1}^n)$. Indeed, there are exactly $k$ particles arriving in the interval $[0, t]$, with the $k$-th one arriving in $[s, t]$, and the remaining $n - k$ particles surviving beyond $t$. Using the survival probability $S(t) = \mathbb{P}(\tau > t)$ of a single Brownian motion in the domain, and the arrival time $\tau^{(i)}$  of the $i$-th particle,  we have
{\small
\begin{align*}
    Pr\lbrace\tau_k^n \in [s,t], \tau_{k+1}^n \geq t\rbrace
    & = \sum_{i_1,...,i_k,  \ l_0 \in[1,k]}Pr\lbrace\tau^{(i_l)}<\tau^{(i_{l_0})}, \tau^{(i_{l_0})} \in [s,t], \forall j \not \in I, \tau^{(j)}>t\rbrace \\
    & = \sum_{I, \  l_0} Pr\lbrace \tau > t\rbrace^{n-k} \int_s^t Pr\lbrace\tau<u\rbrace^{k-1}Pr\lbrace\tau=u\rbrace du \\
    &= \binom{n}{k}Pr\lbrace \tau > t\rbrace^{n-k}(Pr\lbrace\tau<t\rbrace^k-Pr\lbrace\tau<s\rbrace^k) \\
    &= \binom{n}{k}S(t)^{n-k}((1-S(t))^k-(1-S(s))^k),
\end{align*}}
where $\tau^{(i)}$ is the arrival time of the Brownian of index i, $S(s) =  Pr\lbrace \tau>s\rbrace$ is the survival probability of a single Brownian motion in the domain.
Thus, for $s\leq t$, differentiating with respect to $s$ and $t$ gives the joint density of $(\tau_k^n, \tau_{k+1}^n)$ for $s \leq t$,
\begin{equation} \label{order}
    f_{\tau_k^n, \tau_{k+1}^n}(s, t)
    = k(n-k) \binom{n}{k} \, f(s) \, f(t) \, (1 - S(s))^{k-1} S(t)^{n-k-1},
\end{equation}
where $f(t) = -S'(t)$ is the first-passage time density for a single particle.
\begin{figure}[http!]
\begin{adjustwidth}{-1cm}{-1cm}
\centering
\includegraphics[scale = 0.5]{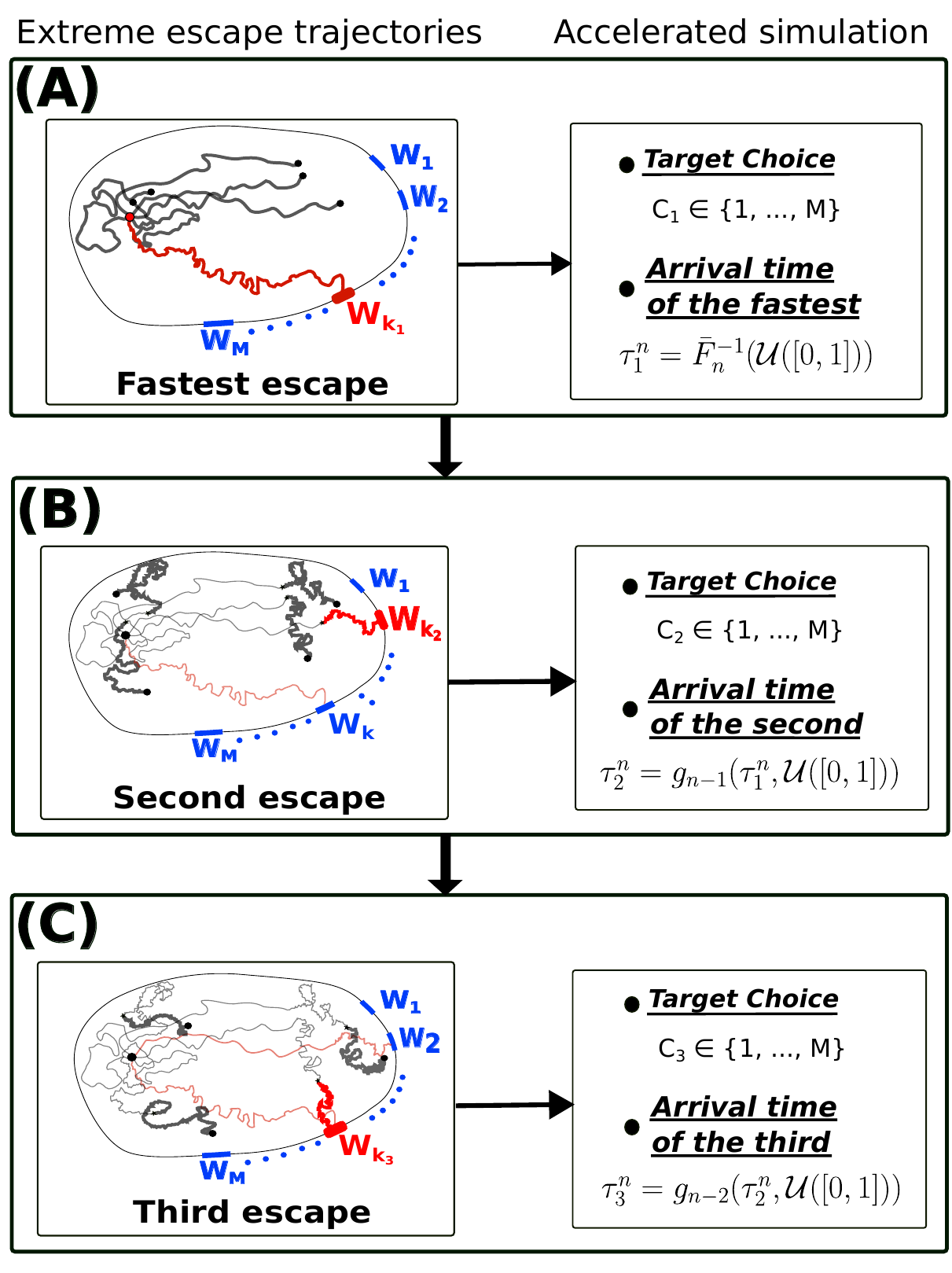}
\caption{\textbf{Iterative Sampling of Successive Escape Events}. Recursive procedure used to simulate the first $k = 3$ escape events: given the arrival time $\tau_k$ of the $k$-th escape, $\tau_{k+1}$ is sampled conditionally, reflecting the decreasing number of particles in the system.}
\end{adjustwidth}
\end{figure}
\subsection{Asymptotic Approximation of successive arrivals}
\begin{figure}[http!]
\centering
\includegraphics[scale = 0.5]{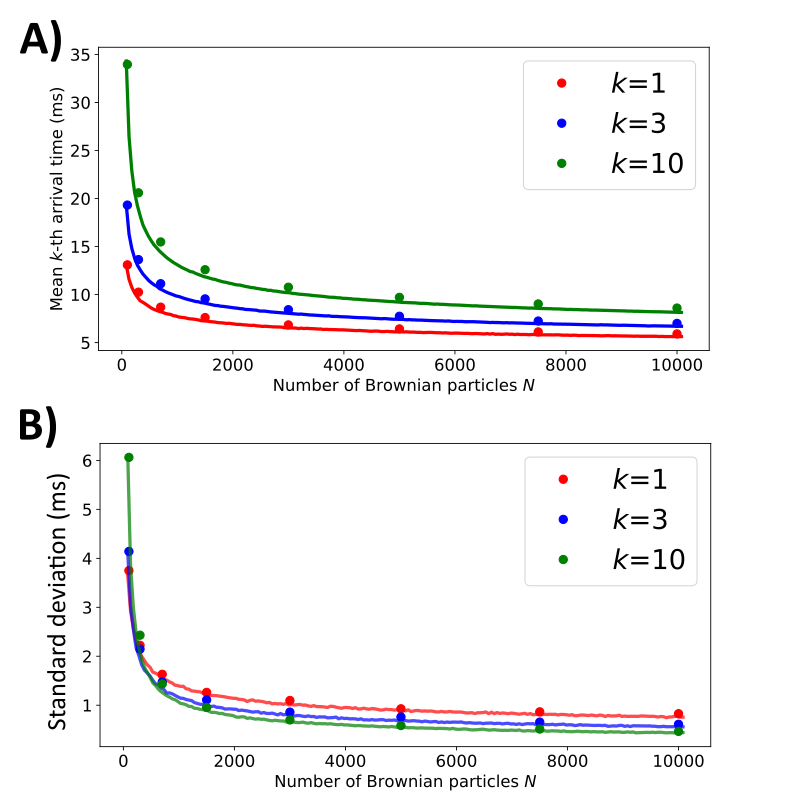}
\caption{\textbf{Arrival Times of the $k$-th Particle : Stochastic Simulation vs Accelerated Algorithm realizations.} The markers are averages over simulations and the lines are average over algorithm realizations for the mean (\textbf{A}) and the variance (\textbf{B}) in one dimension.}\label{figkparticles}
\end{figure}

To evaluate $F^{-1}$, we again use short-time asymptotic approximations (eq. \ref{asympexpressions})for the survival probability $S(t)$ in dimensions 1, 2, and 3. With $S(t)^n \approx \exp(-nF(t))$, we shall derive here the successive arrival time expression:
\begin{equation}\label{tkrecursion}
    t_{k+1} =
    \begin{cases}
      \ds \frac{\delta^2}{2D\, W_0 \left( \ds\frac{2(n-k)^2}{\pi \left[ (n-k)F(t_k) + \log(1/u) \right]^2} \right)} & \text{(1D)} \\\\
        \ds \frac{\delta^2}{4D\, W_0 \left(\ds \frac{\sqrt{2}\pi(n-k)}{8 \log(1/\varepsilon) \left[ \log(1/u) + (n-k)F(t_k) \right]} \right)} & \text{(2D)} \\\\
        \ds -\frac{\delta^2}{2D\, W_{-1} \left(\ds -\frac{\pi \delta^4 \left[ \log(1/u) + (n-k)F(t_k) \right]^2}{2(n-k)^2 a^4} \right)} & \text{(3D)}
    \end{cases}
\end{equation}
where $W_0$ and $W_{-1}$ are the principal and lower branches of the Lambert $W$ function.\\
To  derive explicit inversion formulas for the arrival time \(t_{k+1}\) based on the short-time asymptotics of the cumulative escape probability \(F(t)\), we start from the recursion eq. \ref{recursion}
so that the \((k+1)\)-st arrival time \(t_{k+1}\) is determined by inverting \(F\) at the updated value on the right-hand side.
Using \(F(t)=1-S(t)\), the exit probability \(F(t)=1-S(t)\) is given by:
\begin{equation}
    F(t) \approx
    \begin{cases}
        \frac{\sqrt{4Dt}}{\delta \sqrt{\pi}}\exp\left(-\frac{\delta^2}{4Dt}\right) & \text{(1D)}, \\[0.7em]
        \frac{\sqrt{2}\pi Dt}{2\log\left(\frac{1}{\varepsilon}\right)\delta^2}\exp\left(-\frac{\delta^2}{4Dt}\right) & \text{(2D)}, \\[0.7em]
        \frac{a^2}{\delta\sqrt{\pi Dt}}\exp\left(-\frac{\delta^2}{4Dt}\right) & \text{(3D)}.
    \end{cases}
\end{equation}
\label{eq:exitprob}
We only derive now the inversion for the 1D case. Starting from
\[
F(t)=\frac{\sqrt{4Dt}}{\delta\sqrt{\pi}}\exp\!\left(-\frac{\delta^2}{4Dt}\right),
\]
The short-time asymptotic for the cumulative distribution of the first-passage time in dimension 1 is
\begin{equation}
    F(t) \approx \frac{1}{\sqrt{\pi y}} e^{-y}, \qquad y := \frac{\delta^2}{4Dt}.
\end{equation}
From the recursive formula for the ordered arrival times:
\begin{equation}
    \Delta_{k+1} := F(t_{k+1}) = F(t_k) + \frac{\log\left( \frac{1}{U_{k+1}} \right)}{n - k},
\end{equation}
we substitute the asymptotic form of \( F(t_{k+1}) \) to obtain:
\begin{equation}
    \Delta_{k+1} = \frac{1}{\sqrt{\pi y}} e^{-y}.
\end{equation}
We isolate the exponential term:
\begin{equation}
    \Delta_{k+1} \sqrt{\pi y} = e^{-y} \quad \Rightarrow \quad y e^y = \frac{1}{\pi \Delta_{k+1}^2}.
\end{equation}
This yields the Lambert equation, whose solution is:
\begin{equation}
    y = W_0\left( \frac{1}{\pi \Delta_{k+1}^2} \right),
\end{equation}
where \( W_0 \) is the principal branch of the Lambert function. Substituting back for \( y \) and solving for \( t_{k+1} \), we obtain:
\begin{equation}
    t_{k+1} = \frac{\delta^2}{4D \cdot W_0\left( \dfrac{1}{\pi \left( F(t_k) + \dfrac{\log(1/U_{k+1})}{n - k} \right)^2 } \right)}.
\end{equation}
Similar derivations allow obtaining the expression \ref{tkrecursion} in dims 2 and 3.
The recursive simulation procedure is summarized in the following algorithm:
\begin{algorithm}[http!]
\caption{Recursive Simulation of the First $k$ Arrival Times}
\footnotesize
\begin{algorithmic}[1]
\Require $n$: total number of particles, $k$: number of arrivals to simulate
\Require $\delta_i$: geodesic distances to $M$ targets
\Ensure \texttt{arrival\_times[1:k]}, \texttt{selected\_targets[1:k]}

\State $F \gets 0$
\For{$i = 1$ to $k$}
    \State Draw $U_{2i} \sim \mathcal{U}[0, 1]$
    \State $F \gets F - \log(U_{2i}) / (n - i + 1)$
    \State $t_i \gets F^{-1}(F)$ using dimension-specific Lambert $W$ inversion
    \State Draw $U_{2i+1} \sim \mathcal{U}[0, 1]$
    \State Select target index $j$ via geometric splitting probabilities
    \State \texttt{arrival\_times[i]} $\gets t_i$, \texttt{selected\_targets[i]} $\gets j$
\EndFor
\State \Return \texttt{arrival\_times}, \texttt{selected\_targets}
\end{algorithmic}
\end{algorithm}
{Details of the algorithm is presented in Appendix A and SI.} To validate the algorithm, we compare the mean $k$-th arrival times predicted by the algorithm against values obtained from direct simulations of Brownian particles. As shown in Fig.\ref{figkparticles}, the recursive sampling approach achieves excellent agreement in dimension one, particle numbers, and ranks $k$. This numerical validation illustrates the accuracy of the accelerated algorithm for capturing sequential escape dynamics of the fastest.
\section{Extreme First-Passage under Uniform Killing Constraints}
In this section, we extend our acceleration algorithm to incorporate a uniform killing, where each Brownian particle is assigned a random lifetime drawn independently from an exponential distribution $\tau^{\text{death}} \sim \mathrm{Exp}(\gamma)$ with parameter $\gamma$. If a particle escape time exceeds its lifetime, the particle is terminated before reaching the target and is excluded from the successful escapes.\\
The $k$-th prospective escape time $\tau'_k$ is generated from the recursion (without killing), and the corresponding exponentially distributed lifetime is $\tau_k^{\mathrm{death}} \sim \mathrm{Exp}(\gamma)$ for that particle: if $\tau'_k < \tau_k^{\mathrm{death}}$, the particle has escaped,  otherwise, the particle is discarded, and the remaining $(n - k)$ particles continue to be processed. The arrival time generation still follows the same recursive rule as in the no-killing case, that is
\begin{align*}
\ds F(\tau'_{k+1}) &= F(\tau'_k) + \frac{\log\left(\frac{1}{U_{k+1}}\right)}{n - k}, \\
\tau'_{k+1} &= F^{-1}\left(F(\tau'_k) + \frac{\log\left(\frac{1}{U_{k+1}}\right)}{n - k} \right)= g_{n-k}(\tau'_k, U_{k+1}),
\end{align*}
where $U_{k+1} \sim \mathcal{U}[0,1]$. For each such candidate $\tau'_{k+1}$, we generate an independent exponential lifetime:
\begin{equation}
\tau^{\mathrm{death}}_{k+1} \sim \mathrm{Exp}(\gamma).
\end{equation}
Only candidates satisfying $\tau'_{k+1} < \tau^{\mathrm{death}}_{k+1}$ are retained and the $(k+1)$-th actual escape time is then determined as:
\beq
    \tau_{k+1} = \min \lbrace \tau'_i  \  | \ \tau'_i > \tau_k, \ \tau'_i < \tau^{death}_i \rbrace.
\eeq
To derive the distribution of observed arrival times under killing, we condition on the candidate escape time being less than the corresponding lifetime. The survival pdf of the killing process is:
\begin{equation}
\mathbb{P}(\tau < \tau^{\mathrm{death}}) = \int_0^\infty f_{\tau'}(t) \cdot \mathbb{P}(\tau^{\mathrm{death}} > t) dt = \int_0^\infty f_{\tau'}(t) e^{-\gamma t} dt.
\end{equation}
This gives the modified (effective) density of observed arrival times:
\begin{equation}
f_{\tau}(t) = f_{\tau'}(t) \cdot e^{-\gamma t}.
\end{equation}
This exponential re-weighting penalizes longer escape times and favors earlier arrivals, effectively sharpening the arrival time distribution. The recursive arrival generation proceeds identically, but only those times that satisfy $\tau_i < \tau^{\mathrm{death}}_i$). Finally, once $k$ valid escape times $\tau_1, \tau_2, \dots, \tau_k$ have been collected, they are sorted to construct the proper order of observed arrivals.
This extension preserves the recursive structure of the original algorithm, while incorporating termination. To conclude, the simulation algorithm proceeds as follows:
\begin{enumerate}
    \item Initialize an empty list of successful escape times $\mathcal{T} = \emptyset$.

    \item \textbf{Repeat} until $|\mathcal{T}| = k$:
    \begin{enumerate}
        \item Sample an arrival time $\tau' \sim f_{\tau'}(t)$ from the first-passage time distribution of a Brownian particle in the domain.

        \item Sample an independent killing time $\tau^{\mathrm{death}} \sim \mathrm{Exp}(\gamma)$.

        \item \textbf{If} $\tau' < \tau^{\mathrm{death}}$, \textbf{then} add $\tau'$ to $\mathcal{T}$.
    \end{enumerate}

    \item Sort the accepted arrival times: $\mathcal{T} = \{\tau_1, \tau_2, \dots, \tau_k\}$ with $\tau_1 < \tau_2 < \cdots < \tau_k$.
\end{enumerate}
\begin{figure}[http!]
\centering
\includegraphics[scale = 0.7]{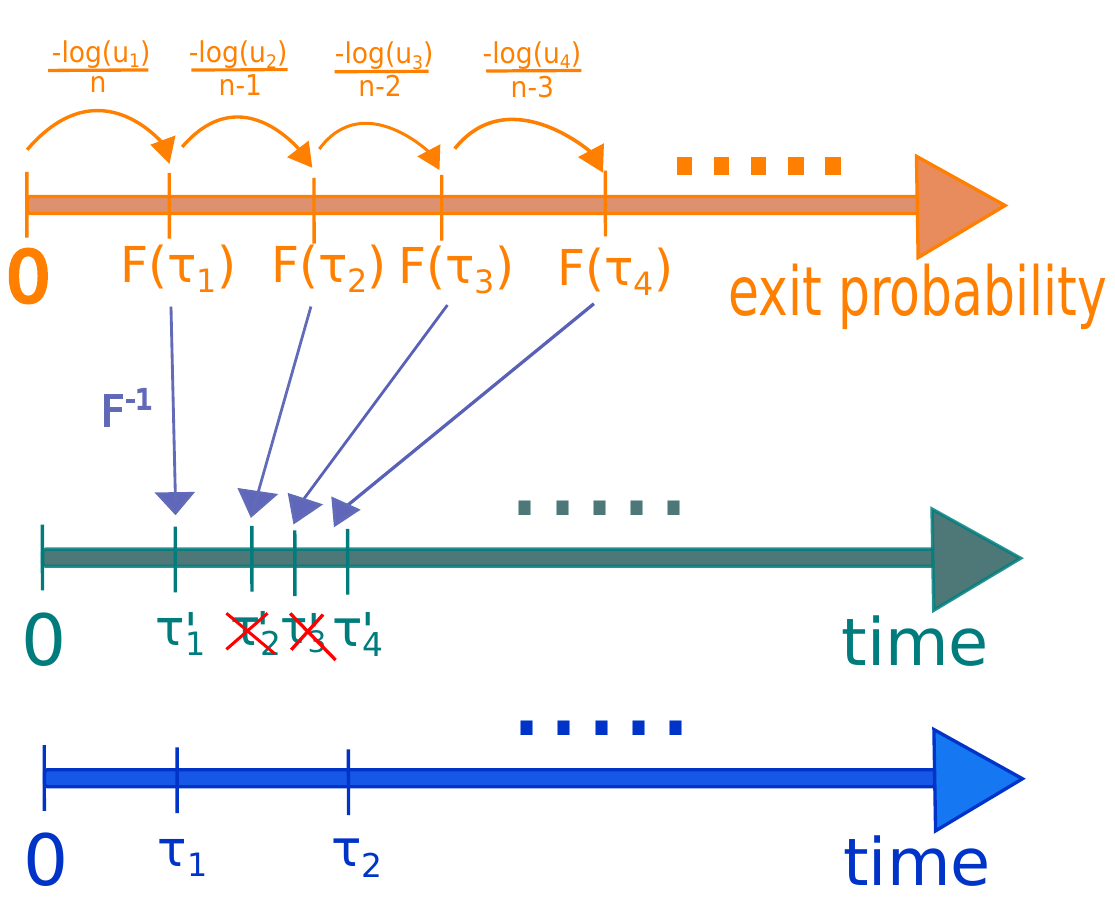}
\caption{\textbf{Arrival Time Sampling under Uniform Killing}. Each Brownian particle is assigned an exponential lifetime. Candidate escape times generated by the recursive inverse transform method are rejected if they exceed the corresponding killing time. Surviving events are retained and sorted by arrival.}
\end{figure}
\section{Asymptotic Behavior of the Mean Fastest Arrival Time under Time-Dependent Emission}
\subsection{Model Setup and Definitions}
We consider a diffusion process on the half-line \( \mathbb{R}^+ = [0, \infty) \), where particles are emitted from position \( x = y > 0 \) and absorbed at the origin \( x = 0 \). The particles diffuse with diffusion constant \( D \), and the emission is governed by a time-dependent rate \( \phi(t) \).  { The emission profile $\phi(t)$ is a nonnegative integrable function on $[0,\infty)$ such that the convolution with the first-passage distribution is well-defined. We shall use here for the emission profile a  gamma  distribution :}
\beq \label{injection}
    \phi(t) = \alpha^2 t e^{-\alpha t}, \quad \alpha > 0,
\eeq
which models a sharp peak centered around \( t = 1/\alpha \). The survival probability of a single particle initially at \( x = y \), diffusing in \( \mathbb{R}^+ \) with absorption at the origin, is given by \cite{Holcmanschuss2018}:
\beq
S(t) = 1 - \frac{\sqrt{4Dt}}{y\sqrt{\pi}} \exp\left(-\frac{y^2}{4Dt}\right).
\eeq
For such stochastic emission following  \( \phi(t) \), the effective survival probability is given by:
\beq
    S_\phi(t) = 1 - \int_0^t \left[1 - S(t - s)\right] \phi(s) \, ds.
\eeq
The MFAT among \( n \) identically emitted and independent particles is then:
\begin{equation} \label{emission}
    \bar{\tau}_{\phi,n} = \int_0^\infty \left[ S_\phi(t) \right]^n dt \approx \int_0^\infty \exp\left( -n [1 - S_\phi(t)] \right) dt.
\end{equation}
We shall present an algorithm to simulate the arrival of the fastest based on the asymptotic behavior of  fastest \( \bar{\tau}_{\phi,n} \) in various limits, starting with the one of \textbf{slow emission}, i.e., \( \alpha \ll \frac{1}{\tau_d \sqrt{n}} \), where the diffusive timescale is:
\beq
\tau_d = \frac{y^2}{4D}.
\eeq
\subsection{Asymptotic Expansion in the Slow Emission Regime}
Assuming \( \bar{\tau}_{\phi,n} \gg \tau_d \), we approximate \( e^{-\alpha t} \approx 1 \) in the expression of \( \phi(t) \). We first expand \( 1 - S_\phi(t) \) using the convolution:
\begin{align}
    1 - S_\phi(t) &\approx \int_0^t \alpha^2 s \cdot \frac{\sqrt{4D(t - s)}}{y \sqrt{\pi}} \exp\left(-\frac{y^2}{4D(t - s)}\right) ds \\
    &= \int_0^t \alpha^2 s \exp\left(-\frac{\tau_d}{t - s}\right) \sqrt{\frac{t - s}{\pi \tau_d}} \, ds.
\end{align}
With the substitution \( s = tu \), this becomes:
\begin{align}
    1 - S_\phi(t) &\approx \frac{\alpha^2 t^{5/2}}{\sqrt{\pi \tau_d}} \int_0^1 u \sqrt{1 - u} \exp\left(-\frac{\tau_d}{t(1 - u)}\right) du.
\end{align}
In the limit \( t \gg \tau_d \), the exponential factor is sharply peaked near \( u = 0 \), and we approximate:
\beq
    1 - S_\phi(t) \approx \frac{\alpha^2 t^{5/2}}{\sqrt{\pi \tau_d}} \exp\left( -\frac{\tau_d}{t} \right) \int_0^1 u \sqrt{1 - u} \, du.
\eeq
Evaluating the integral:
\[
\int_0^1 u \sqrt{1 - u} \, du = \frac{4}{15},
\]
we obtain:
\beq
    1 - S_\phi(t) \approx \frac{4\alpha^2 \tau_d^2}{15\sqrt{\pi}} \left( \frac{t}{\tau_d} \right)^{5/2} \exp\left( -\frac{\tau_d}{t} \right).
\eeq
Thus, the MFAT becomes:
\begin{align}
    \bar{\tau}_{\phi,n} &\approx \int_0^\infty \exp\left( - \frac{4n \alpha^2 \tau_d^2}{15\sqrt{\pi}} \left( \frac{t}{\tau_d} \right)^{5/2} \exp\left( -\frac{\tau_d}{t} \right) \right) dt.
\end{align}
Introducing the change of variable \( u = t / \tau_d \), and defining:
\beq
    \bar{\alpha} := \frac{4n \alpha^2 \tau_d^2}{ 15\sqrt{\pi}},
\eeq
we rewrite the integral as:
\beq
    \bar{\tau}_{\phi,n} = \tau_d \int_0^\infty \exp\left( - \bar{\alpha}^2 u^{5/2} e^{-1/u} \right) du.
\eeq
Performing a scaling \( u = v / \bar{\alpha}^{4/5} \), we obtain:
\begin{align}
    \bar{\tau}_{\phi,n} &= \frac{\tau_d}{\bar{\alpha}^{4/5}} \int_0^\infty \exp\left( - v^{5/2} \exp\left( -\frac{\bar{\alpha}^{4/5}}{v} \right) \right) dv.
\end{align}
In the limit \( \bar{\alpha} \to 0 \), the exponential in the integrand converges to unity, yielding:
\[
\int_0^\infty \exp(-v^{5/2}) dv = \frac{2}{5} \Gamma\left( \frac{2}{5} \right),
\]
so that:
\beq
    \bar{\tau}_{\phi,n} \sim \frac{2\tau_d}{5\bar{\alpha}^{4/5}} \Gamma\left( \frac{2}{5} \right).
\eeq
Substituting \( \bar{\alpha} \), we arrive at the final asymptotic expression in terms of the original parameters:
\beq
    \bar{\tau}_{\phi,n} \sim \left( \frac{9 \pi}{250} \right)^{1/5} \cdot \frac{y^{2/5}}{D^{1/5} \alpha^{4/5} n^{2/5}} \cdot \Gamma\left( \frac{2}{5} \right)
\eeq
This asymptotic formula is valid in the regime \( \alpha \ll \frac{1}{\tau_d \sqrt{n}} \), or equivalently:
\[
\alpha \ll \frac{4D}{y^2 \sqrt{n}},
\]
and implies that:
$ \bar{\tau}_{\phi,n} \gg \tau_d = \frac{y^2}{4D}.$
\begin{figure}[http!]
    \begin{adjustwidth}{-2cm}{-1cm}
    \centering
    \includegraphics[scale=0.9]{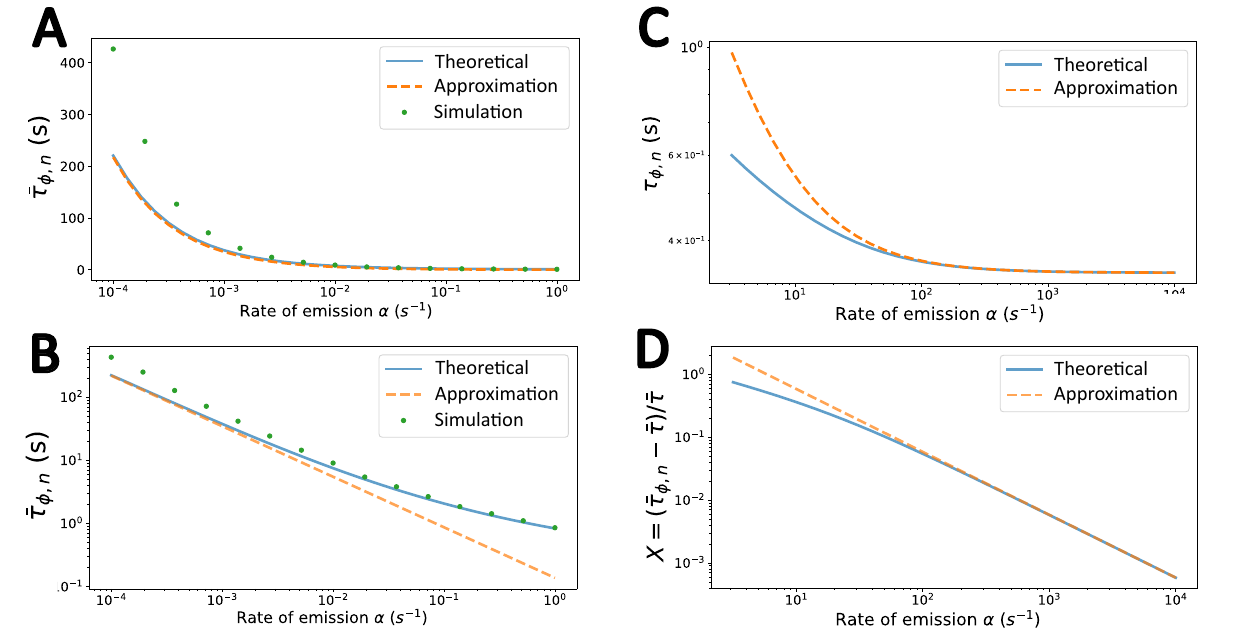}
    \end{adjustwidth}
    \begin{adjustwidth}{-1cm}{-1cm}
    \caption{\textbf{Asymptotic behavior of the MFAT in the slow and fast emission regimes vs accelerated simulations.} The slow regime is plotted for $N=1000$ in \textbf{A} and \textbf{B} in different scales. The asymptotic approximation (eq. \ref{eq:exitprob} dotted curve) agrees well with theoretical formula (eq. \ref{distrbutiontaun}) for small \(\alpha\) limit. Discrepancies appear as \(\alpha \to 0\), due to the approximation of short-time asymptotics for the expression of \(S(t)\) beyond their validity range. The fast regime is plotted for $N=1000$ in \textbf{C} and \textbf{D} where the relative difference between the fast emission regime and the instantaneous one is shown.}
    \end{adjustwidth}
    \label{fig5}
\end{figure}
{ The non-monotonic dependence of the MFAT on the emission rate $\alpha$ (Fig. \ref{fig5}) can be understood as follows: in the slow emission regime ($\alpha \ll 1/\tau_d$), particles are injected over a timescale much longer than the diffusive escape time, so the fastest arrival is delayed by the emission process itself. Conversely, in the fast emission regime $\alpha \gg 1/\tau_d$, all particles are released nearly instantaneously, and the MFAT converges to the instantaneous-emission limit. The intermediate regime reflects a competition between these two timescales.}
\subsection{Asymptotic Regimes: Fast and Intermediate Emission Profiles}
\subsubsection*{Fast Emission Regime: \(\alpha \gg \frac{\log(n)^2}{\tau_d}\)}
We now consider the opposite limit of the slow emission case, where the emission (eq. \ref{injection}) occurs on a timescale much faster than diffusion. To compute the mean time and obtain the fast regime inequality, we start from \( f(t) = 1 - S(t) \) and
\beq
    f(t) = \frac{y}{\sqrt{\pi}} \sqrt{ \frac{4D}{t} } \exp\left( -\frac{y^2}{4Dt} \right).
\eeq
Then, the effective survival correction is given by:
\begin{align}
    1 - S_\phi(t) &= \int_0^t e^{-\alpha s} \alpha^2 s f(t - s) \, ds.
\end{align}
Expanding \( f(t - s) \) around \( s = 0 \), we use the Taylor expansion:
\beq
    f(t - s) \approx f(t) - s f'(t),
\eeq
to obtain:
\begin{align}
    1 - S_\phi(t) &\approx f(t) \int_0^t e^{-\alpha s} \alpha^2 s \, ds - f'(t) \int_0^t e^{-\alpha s} \alpha^2 s^2 \, ds \\
    &\approx f(t) \int_0^\infty e^{-\alpha s} \alpha^2 s \, ds - f'(t) \int_0^\infty e^{-\alpha s} \alpha^2 s^2 \, ds \\
    &= f(t) - \frac{2}{\alpha} f'(t),
\end{align}
where we used
\[
\int_0^\infty \alpha^2 s e^{-\alpha s} \, ds = 1, \qquad \int_0^\infty \alpha^2 s^2 e^{-\alpha s} \, ds = \frac{2}{\alpha}.
\]
At short times \( t \ll \tau_d \), we use the approximation \( f'(t) \approx \frac{\tau_d}{t^2} f(t) \), yielding:
\beq \label{eq:effective_survival_fast}
    1 - S_\phi(t) \approx f(t) \left( 1 - \frac{2 \tau_d}{\alpha t^2} \right).
\eeq
Therefore, the MFAT becomes:
\beq \label{eq:mfat_integral}
    \bar{\tau}_{\phi,n} \approx \int_0^\infty \exp\left( -n f(t) \left( 1 - \frac{2 \tau_d}{\alpha t^2} \right) \right) dt.
\eeq
Following the approach in \cite{grebenkov2025}, we approximate the integrand by a Heaviside function. The MFAT is then determined by the location \(\bar{\tau}_{\phi,n}\) such that the argument of the exponential reaches a fixed threshold \( \beta \in (0,1) \). Indeed, for large $n$, the exponential term in Eq.~\eqref{eq:mfat_integral} becomes sharply peaked around the time where its argument becomes of order one. Following the asymptotic approach  \cite{grebenkov2025}, we approximate the integrand by a Heaviside step function:
\begin{equation}
\exp\left( -n f(t) \left(1 - \frac{2 \tau_d}{\alpha t^2} \right) \right) \approx \Theta(\bar{\tau}_{\phi,n} - t),
\label{eq:heaviside_approx}
\end{equation}
where $\Theta$ is the Heaviside function, and $\bar{\tau}_{\phi,n}$ is the time where the exponential argument reaches a small fixed threshold $\beta \in (0,1)$:
\begin{equation}
n f(\bar{\tau}_{\phi,n}) \left(1 - \frac{2 \tau_d}{\alpha \bar{\tau}_{\phi,n}^2} \right) = \beta.
\label{eq:threshold_equation}
\end{equation}
This condition defines the point where the integrand transitions from values near 1 to near 0, making it the dominant contribution to the integral in Eq.~\eqref{eq:mfat_integral}. Solving Eq.~\eqref{eq:threshold_equation} gives an asymptotic estimate of the MFAT. To match the unperturbed case (i.e., with instantaneous emission $\phi(t) = \delta(t)$), where  $\bar{\tau}$ corresponds MFAT:
\begin{equation}
\bar{\tau} := \text{MFAT for instantaneous emission}.
\end{equation}
Matching the two sides gives:
\beq
    n f(\bar{\tau}_{\phi,n}) \left( 1 - \frac{2 \tau_d}{\alpha \bar{\tau}_{\phi,n}^2} \right) = n f(\bar{\tau}),
\eeq
where \( \bar{\tau} \) is the MFAT in the case of instantaneous emission (i.e., \(\phi(t) = \delta(t)\)). Since \( \bar{\tau}_{\phi,n} \to \bar{\tau} \) as \( \alpha \to \infty \), we perform a first-order Taylor expansion:
\[
\bar{\tau}_{\phi,n} = \bar{\tau}(1 + X), \qquad \text{with} \quad X \ll 1.
\]
Using \( f(\bar{\tau}_{\phi,n}) \approx f(\bar{\tau})(1 + X \frac{\tau_d}{\bar{\tau}}) \) and inserting into the identity, we get:
\begin{align}
    f(\bar{\tau})\left(1 + X \frac{\tau_d}{\bar{\tau}}\right) \left( 1 - \frac{2\tau_d}{\alpha \bar{\tau}^2} \right) &\approx f(\bar{\tau}), \\
    1 + X \frac{\tau_d}{\bar{\tau}} &\approx 1 + \frac{2\tau_d}{\alpha \bar{\tau}^2}, \\
    X &\approx \frac{2}{\alpha \bar{\tau}}.
\end{align}
Hence, the corrected MFAT satisfies:
\beq
    \bar{\tau}_{\phi,n} \approx \bar{\tau} + \frac{2}{\alpha}
\eeq
This requires the emission to be significantly faster than the typical diffusion-driven escape time.  From Eq.~\eqref{eq:effective_survival_fast}, we see that emission smooths the peak by a factor of $1 - \frac{2 \tau_d}{\alpha t^2}$. To retain sharp decay, we require this correction to be small at $t = \bar{\tau}$:
\begin{equation}
\frac{2 \tau_d}{\alpha \bar{\tau}^2} \ll 1 \quad \Rightarrow \quad \alpha \gg \frac{2 \tau_d}{\bar{\tau}^2}.
\end{equation}
Using that $\bar{\tau} \sim \tau_d / \log(n)$ in the fast arrival limit, we obtain the final condition:
\begin{equation}
\alpha \gg \frac{\log(n)^2}{\tau_d},
\label{eq:alpha_condition}
\end{equation}
which ensures that the emission kernel $\phi(t)$ acts nearly as an impulse compared to the diffusion timescale, thereby validating the fast emission regime.
\subsubsection*{Intermediate Regime}
To interpolate between the slow and fast emission regimes, we consider the effective survival probability derived in the slow emission limit and define \( \bar{\tau}_{\phi,n} \) by setting:
\beq
    1 - S_\phi(\bar{\tau}_{\phi,n}) = \frac{\log(2)}{n}.
\eeq
Using the slow emission asymptotic form:
\beq
    1 - S_\phi(t) \approx \frac{4 \alpha^2 \tau_d^2}{15\sqrt{\pi}} \left( \frac{t}{\tau_d} \right)^{5/2} \exp\left( -\frac{\tau_d}{t} \right),
\eeq
we obtain:
\begin{align}
    \frac{4 n \alpha^2 \tau_d^2}{15\sqrt{\pi}} \left( \frac{\bar{\tau}_{\phi,n}}{\tau_d} \right)^{5/2} \exp\left( -\frac{\tau_d}{\bar{\tau}_{\phi,n}} \right) = \log(2).
\end{align}
Using \( z = \frac{2\tau_d}{5 \bar{\tau}_{\phi,n}} \). Then the equation becomes:
\begin{align}
    z e^z = \frac{2}{5} \left( \frac{\bar{\alpha}^2}{\log(2)} \right)^{2/5},
\end{align}
where \( \bar{\alpha}^2 = \frac{4n \alpha^2 \tau_d^2}{15\sqrt{\pi}} \). Solving using the principal branch \( W_0 \) of the Lambert W function, we find:
\beq
    \bar{\tau}_{\phi,n} \approx \frac{2\tau_d}{5 W_0\left( \frac{2}{5} \left( \frac{4n\alpha^2 \tau_d^2}{15\sqrt{\pi} \log(2)} \right)^{2/5} \right) },
\eeq
in the orginal variable,
\beq
\bar{\tau}_{\phi,n} \approx \frac{y^2}{10D} \cdot \frac{1}{W_0 \left( \frac{2}{5} \left( \frac{n \alpha^2 y^4}{60 D^2 \sqrt{\pi} \log(2)} \right)^{2/5} \right)}.
\eeq
To conclude, this expression provides a smooth interpolation between the slow and fast emission regimes.
\subsection{Algorithm for Simulating Extreme Statistics with Time-Dependent Injection}
To extend the simulation of extreme statistics to the case of non-instantaneous, time-dependent emission, we adapt the recursive sampling method previously developed for instantaneous release. In this setting, the emission profile \( \phi(t) \) defines the injection of particles over time, thereby modifying the effective survival probability. To simulate the successive order statistics of arrival times, we replace the standard cumulative distribution function \( F(t) = 1 - S(t) \) by incorporating emission: the effective survival function \( S_\phi(t) \). Accordingly, the iterative sampling of the \( (k+1) \)-th arrival time \( t_{k+1} \) is governed by the following nonlinear integral equation:
\beq
    \int_{0}^{t_{k+1}} F(t_{k+1} - s)\, \phi(s)\, ds = \int_{0}^{t_k} F(t_k - s)\, \phi(s)\, ds - \frac{\log\left(\frac{1}{U_{k+1}}\right)}{n - k},
\eeq
where \( U_{k+1} \sim \mathcal{U}(0,1) \) is an independent uniform random variable. This relation arises from conditioning on the cumulative escape probability of the remaining \( n-k \) particles after the \( k \)-th arrival. The left-hand side reflects the effective probability that a particle escapes before time \( t_{k+1} \), convolved with the emission kernel \( \phi(s) \), while the right-hand side tracks the cumulative contribution from the previous \( k \) escapes and adjusts it via a random exponential increment, consistent with order statistics theory. To compute \( t_{k+1} \), a numerical root-finding algorithm (such as  Newton's method) can be used. This framework allows for accurate sampling of ordered arrival times in systems with temporally structured injections and offers an efficient alternative to direct simulation of Brownian trajectories under emission dynamics. The associated pseudocode is given by:
\begin{algorithm}[H]
\caption{Recursive Simulation with Time-Dependent Emission}
\footnotesize
\begin{algorithmic}[1]
\Require $n$: number of particles, $k$: number of arrivals to simulate
\Require $\phi(t)$: emission profile (e.g., bell-shape or gamma-distribution)
\Require $F(t)$: cumulative distribution of first-passage time
\Ensure \texttt{arrival\_times[1:k]}

\State Initialize: $t_0 \gets 0$, $C_0 \gets 0$ \Comment{Cumulative integral term}
\For{$i = 1$ to $k$}
    \State Draw $U_i \sim \mathcal{U}(0,1)$
    \State Set target cumulative: $C_i = C_{i-1} - \frac{\log(1/U_i)}{n - i + 1}$
    \State Define function $G(t)$:
    \[
        G(t) = \int_0^t F(t - s)\phi(s)\, ds - C_i
    \]
    \State Solve for $t_i$ such that $G(t_i) = 0$ using root-finding
    \State Store $t_i$ in \texttt{arrival\_times[i]}
\EndFor
\State \Return \texttt{arrival\_times}
\end{algorithmic}
\end{algorithm}
\section{Algorithmic Validation and Evaluation}
We now assess the efficiency of our algorithms, focusing on accuracy and convergence towards the theoretical distributions of arrival times, computational cost and scalability, and robustness under application-driven scenarios.
\subsection{Accelerated simulation in the pure diffusion case}
We use the 1D-diffusion setting to assess the ability of our algorithm to produce the distribution of successive arrival times, comparing our algorithm with the theoretical predictions for different values of $N$ (Fig. ~\ref{fig:algo3}). The error plots illustrates convergence of the algorithm as $N$ increases. The convergence is especially clear for early arrival orders (small $k$ relative to $N$), whereas higher-order arrivals exhibit slower convergence since the approximation $k \ll N$ does not initially hold. Nonetheless, the method remains numerically stable even for very large populations (beyond $N = 10^8$), far beyond what can be simulated directly.\\
{The principal advantage of our algorithm lies in its efficiency: each data point corresponds to $5000$ independent runs of the algorithm for three arrival orders across different values of $N$,  with $5 \times 10^5$ realization in total (the computation completes under 20 seconds on a standard desktop).  In contrast, simulating $5000$ Brownian trajectories (Euler's scheme) trajectories for $N=10^5$ particles with a timestep $10^{-4}$ and average trajectory length $6 \times 10^{-5}$) would require over $3 \times 10^{10}$ steps, each involving Gaussian sampling, boundary checks, and boundary reflections. Such a simulation would take several minutes per data point (Fig. \ref{fig:algo3}). Thus our method achieves orders-of-magnitude savings in computational cost.  \\
In terms of complexity, runtime is essentially independent of $N$, the total number of Brownian particles, but grows linearly in $k$, the arrival order up to which statistics are required. This scaling property makes the algorithm well-suited for large-population problems where only the first few arrivals are of interest.}
\subsection{Applications and Robustness}
To evaluate the applicability and the  robustness of our algorithms, we used the two modified diffusion scenarios: 1- with killing and 2-Time-dependent emission scenario.
\begin{figure}[http!]
    \centering
    \begin{adjustwidth}{-1cm}{-1cm}
    \includegraphics[scale=0.6]{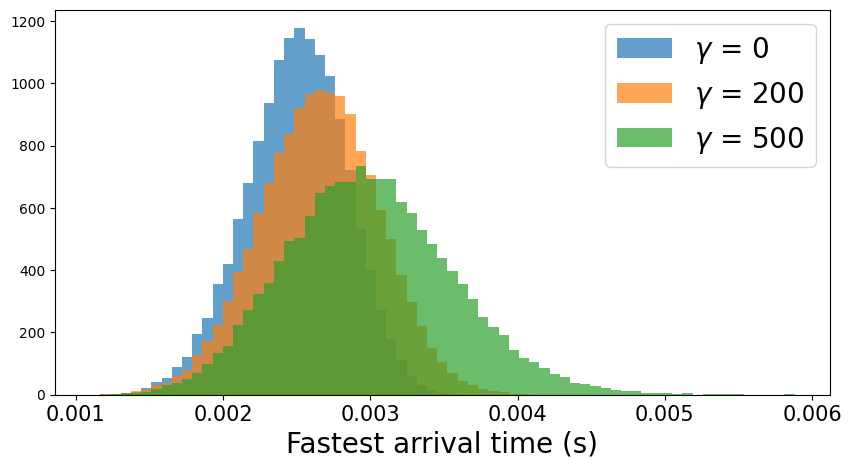}
    \caption{\textbf{Histograms of the first arrival time with killing} with  50000 realizations of our simulation algorithm, each of the first arrival time for three killing rates $\gamma \in \lbrace 0 , 200, 500\rbrace$.  }
    \label{fig:algo2}
    \end{adjustwidth}
\end{figure}
\subsubsection*{Killing scenario}
To study the impact of a killing rate $\gamma$ on the arrival times, we  generated simulations to collect the  histograms of the first arrival time under different values of $\gamma=\{0,200,500\}$ (Fig. \ref{fig:algo2}). The algorithm reliably captures the expected delay in the fastest arrival as $\gamma$ increases, confirming robustness to such perturbations. { The shift of the arrival time distribution with increasing killing rate $\gamma$ has the following interpretation: particles with longer trajectories are more likely to be terminated before reaching the target (Fig. \ref{fig:algo2}). Thus only the fastest diffusers survive, but since some would be early arrivals are now killed, the observed first arrival time increases on average. To conclude,  the killing mechanism acts as a filter that preferentially removes slower trajectories, thereby shifting the observed arrival-time distribution toward faster events and reducing the mean first arrival time.}
\subsubsection*{Time-dependent emission scenario}
\begin{figure}[http!]
    \centering
    \includegraphics[scale=0.6]{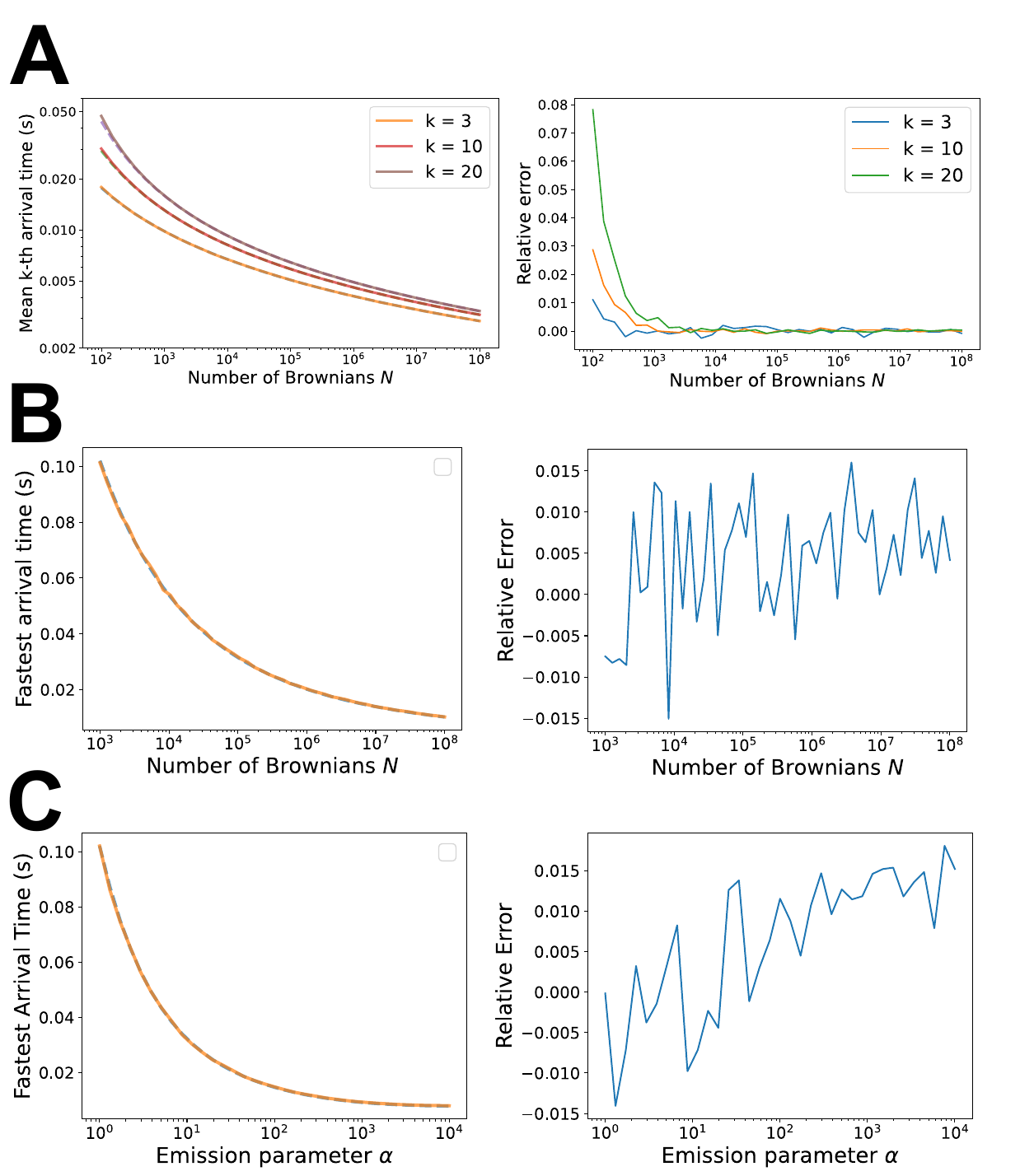}
    \caption{\textbf{Comparison between Theoretical and Simulated Mean Arrival Times under Different Emission Regimes.}
Mean arrival times and their corresponding relative errors between theoretical predictions and accelerated simulations vs number N.
\textbf{(A)} Comparing instantaneous emission algorithm: mean arrival times for the 3\textsuperscript{rd}, 10\textsuperscript{th}, and 20\textsuperscript{th} particles, compared with theoretical values computed from the density of the $k^{\text{th}}$ arrival time (eq.~\ref{order}).
\textbf{(B)} Mean fastest arrival time (MFAT) vs number of particles $N$, averaged over 1000 simulations (theoretical values for time-dependent emission are computed from eq.~\ref{emission}).
\textbf{(C)} MFAT vs emission parameter $\alpha$, averaged over 1000 simulations and compared with theoretical values for time-dependent emission (eq.~\ref{emission}). For each case, the right-hand subplots display the relative errors $\frac{\tau_{\text{sim}}-\tau_{\text{theory}}}{\tau_{\text{theory}}}$ with respect to theoretical predictions.}
    \label{fig:algo3}
\end{figure}

In the emission setting, particles are emitted at a rate governed by the parameter $\alpha$. We compare the asymptotic behavior of the MFAT in the fast emission regime, for different values of $\alpha$ and $N$ (Fig. \ref{fig:algo3}): Despite the additional complexity of inverting the survival probability, the algorithm produces stable results within reasonable runtime (about few minutes for datasets comparable in size to those in Fig. \ref{fig:algo3}).
{  To conclude, the different regimes illustrated in Fig.\ref{fig:algo3} reflect the competition between the emission timescale ($\sim 1/\alpha$) and the diffusive timescale $\tau_d$. When emission is slow, the injection process dominates the arrival time, whereas in the fast emission regime the behavior approaches that of instantaneous release.}
\subsubsection*{Limitations}
{A limitation of the present approach lies in its reliance on asymptotic expressions that are most accurate in the large particle-number regime. Accordingly, our analysis and numerical validation have primarily focused on $N \gtrsim 10^2$, with particularly strong agreement observed for $N \gtrsim 10^3$, where relative errors on mean quantities remain below 1\%. For smaller populations, deviations become more pronounced, especially for higher-order statistics such as the $k$-th arrival time with $k$ not negligible compared to $N$ (e.g., the 20$^{th}$ arrival out of $N \sim 100$ particles in Fig.~\ref{{fig:algo3}}). In such cases, the underlying asymptotic assumptions become less accurate, although the proposed method still provides reasonable approximations.\\
A second limitation concerns the scope limited here to standard Brownian diffusion and does not directly apply to non-diffusive transport mechanisms such as active motion, anomalous diffusion, or L\'evy flights. Extending the algorithm to such processes would require deriving analogous short-time survival probability asymptotics, which may not always be available in closed form.}\\
{An additional limitation arises in scenarios involving time-dependent or structured emission profiles, where the survival probability involves convolution-type expressions. The implementation of our algorithm then requires the inversion of such functions to sample arrival times, which we perform using numerical root-finding procedures. This step necessitates sufficient regularity of the emission profile (e.g., smoothness) to ensure numerical stability and efficiency, and may become computationally demanding for highly irregular or discontinuous inputs. Nevertheless, the method remains orders of magnitude faster than direct trajectory-based simulation in all regimes considered here. Finally, the proposed framework relies on asymptotic approximations of first-passage distributions and is therefore most accurate in the large-$N$ regime, where extreme statistics dominate.}
\section{Competing Targets: Splitting Probability in One Dimension}
We now recall the splitting probability that the first to arrive among $n$ independent Brownian particles escapes through one target before another \cite{Schuss:Book,basnayake2018asymptotic}. The splitting probability  reflects the competition between multiple absorbing boundaries and has been studied in various asymptotic regimes \cite{taflia2007dwell,holcman2015stochastic, lawley2022extremehiting}.
{ The fastest particle among many searchers is more likely to reach the closest target \cite{lawley2022extremehiting,lawley2024hittingfast}, as studied based on large-deviation principles, leading to scaling laws $P \sim N^{1-\beta} (\log N)^\rho$ in the large-$N$ limit. We shall derive explicit and constructive formulas for splitting probabilities  including prefactors and corrections that depend on geometric parameters such as distances and target sizes. In addition, we characterize the transition regime near the symmetry point through a boundary-layer analysis, revealing a transition region of width $O(1/\log n)$. In this context, the splitting probabilities are not only analyzed theoretically but are directly used for sampling the identity of the first-arriving particle. This constructive and computational perspective  enables efficient simulation of extreme first-passage events in complex geometries. Here, we focus on dimension 1, with two absorbing points.}\\
We begin with $n$ Brownian particles are initially located at position $y \in \mathbb{R}^+$, with two absorbing points located at $x_1$ and $x_2$ and $a_1 = |y - x_1|$ and $a_2 = |y - x_2|$ are the distances from the source to each target. The arrival times $\tau_1^n$ and $\tau_2^n$ of the first particle to reach $x_1$ or $x_2$, respectively are used to define the splitting probability $\mathbb{P}(\tau_1^n < \tau_2^n)$. We compute the splitting probability from the arrival time distributions $f_i(t) = Pr\left \lbrace \tau_1^n =t \right \rbrace = -\frac{d}{dt}Pr \left \lbrace \tau_i^n > t \right \rbrace$, for $i = 1,2$ \cite{Basnayake2018},
\beq\label{expression1}
Pr \left \lbrace \tau_i^n >t \right \rbrace = \exp \left \lbrace \frac{-N \sqrt{4Dt}e^{-\frac{a_i^2}{4Dt}}}{a_i\sqrt{\pi}}\right \rbrace.
\eeq
The probability is
\beq
Pr \left \lbrace \tau^n_1 < \tau^n_2 \right \rbrace 
&=& \int_{0}^{\infty}\int_{s}^{\infty}f_2(t)f_1(s)dtds= \int_{0}^{\infty}f_1(s) Pr \left \lbrace \tau^n_2 > s \right \rbrace ds \nonumber \\
&=& \int_{0}^{\infty} \exp \left \lbrace \frac{-n \sqrt{4Ds}}{\sqrt{\pi}}\left(\frac{e^{-\frac{a_2^2}{4Ds}}}{a_2}+\frac{e^{-\frac{a_1^2}{4Ds}}}{a_1} \right)\right \rbrace\frac{n \sqrt{4Ds}}{a_1 \sqrt{\pi}}\left(\frac{1}{2s}+\frac{a_1^2}{4Ds^2}\right)e^{-\frac{a_1^2}{4Ds}} ds. \nonumber
\eeq
The change of variable $u = \frac{4Ds}{a_1^2}$, transform the expression above into
\begin{align*}
Pr \left \lbrace \tau^n_1 < \tau^n_2 \right \rbrace
= \frac{n}{\sqrt{\pi}} \int_{0}^{\infty} \exp \left \lbrace \frac{-n}{\sqrt{\pi}}\left( \sqrt{u}e^{-\frac{1}{u}}+\frac{a_1}{a_2}\sqrt{u}e^{-\frac{a_2^2}{a_1^2u}} \right)\right \rbrace e^{-\frac{1}{u}}\sqrt{u}\left(\frac{1}{2u}+\frac{1}{u^2}\right) du \nonumber
\end{align*}
With the additional change of variable \cite{Basnayake2018} $w(u) = \sqrt{u}e^{-\frac{1}{u}}$ where $w'(u) = w \left(\frac{1}{2u}+\frac{1}{u^2}\right) $, for $u$ small, $\ln w \approx -\frac{1}{u}$, we have
\beq
Pr \left \lbrace \tau^n_1 < \tau^n_2 \right \rbrace
\approx: Sp(\lambda, n)
\eeq
where  $\lambda =\frac{a_1}{a_2}$  and
\begin{align*}
Sp(\lambda, n) 
 =: \frac{1}{\sqrt{\pi}}\int_0^\infty \exp\left(-\frac{1}{\sqrt{\pi}}\left(u+\lambda n^{1-\frac{1}{\lambda^2}} u^{\frac{1}{\lambda^2}}\right)\right) dw
 \end{align*}
We recover for any value of $n$, when $a_1 = a_2$ that $Pr \left \lbrace \tau^1 < \tau^2 \right \rbrace = \frac{1}{2}$. In addition, the splitting probability converges towards a Heaviside-like function with a discontinuity in $\lambda = 1$ (see Appendices):
\beq
\lim_{n \rightarrow \infty} Sp(\lambda, n)=H(\lambda)
\eeq
where
\beq
H(\lambda)=
\begin{dcases}
1 \ & \textbf{ for } \ \lambda < 1, \\
\frac{1}{2} \ & \, \textbf{ } \ \lambda = 1, \\
0 \ & \textbf{ } \, \ \lambda > 1.
\end{dcases}
\eeq
The  asymptotic expansions of $Sp(\lambda,n)$ for large $n$ is summarized in the appendix, leading to
\beq
    Sp(\lambda,n) \approx
    \begin{dcases}
         1-\lambda \ \Gamma(1+\frac{1}{\lambda^2})\left( \frac{\sqrt{\pi}}{n}\right)^{\frac{1}{\lambda^2}-1} \\
         \lambda\Gamma(\lambda^2)\left( \frac{\sqrt{\pi}}{n \lambda} \right)^{\lambda^2-1}
    \end{dcases}
\eeq
When $\lambda \approx 1$ (the discontinuity), we use a boundary layer analysis to resolve the profile. Indeed, using the scaled variable
\beq
\lambda = 1 + \frac{X}{ \log(n)} \text{ where } X \sim O(1),
\eeq
by expanding $Sp(\lambda, n)=\tilde Sp(X,n) $ in powers of $u$ with the approximation that $u^\frac{1}{\lambda^2} \approx u$, we have
\begin{align*}
-\frac{1}{\sqrt{\pi}}\left(u+\lambda n^{1-\frac{1}{\lambda^2}}u^{\frac{1}{\lambda^2}} \right) 
    & \approx -\frac{u}{\sqrt{\pi}}\left(1+\exp(2X) \right)
\end{align*}
Thus
\beq
\tilde Sp(X,n)=\frac{1}{1+\exp(2X)},
\eeq
which leads for $\lambda\approx1, Sp(\lambda, n)= \frac{1}{1+n^{2(\lambda-1)}}$, showing the discontinuity of the limit and a band of width $\sim1/\log(n).$
\section{Discussion}
Estimating the timing of rare activation events driven by the fastest among many stochastic trajectories has become a central computation in cellular biophysics and beyond \cite{holcman2015stochastic,Sokolov2019extreme,lawley2022extremehiting}. Over the past two decades, extreme statistics for Brownian particles has been used to describe such rare events, especially in the context of synaptic transmission \cite{schuss2019redundancy}, intracellular calcium signaling \cite{Basnayake2018} or immune receptor selection, where the trigger depends not on average behavior but on the arrival of the first molecule to a specific site.\\
The mean first-passage time to small absorbing regions on the boundary of a bounded domain \cite{Schuss:Book,holcman2015stochastic,benichou2014first,grebenkov2020single} has provided insights into the geometry-dependent scaling of MFPTs and the role of small target sizes but are inherently centered on average arrival properties of single particles. The extension to multiple particles has brought to light a radically different regime: the escape time of the fastest particle among \(n\) diffusing searchers where asymptotic computations \cite{Basnayake2018,schehr2014exact, lawley2020distribution} have shown that the extreme arrival time decays generically logarithmically with \(n\), underscoring the slow convergence of redundancy in search processes, although extended distribution leads to power law behavior. {Recent studies on fast stochastic search \cite{lawley2024hittingfast} have introduced an alternative viewpoint based on conditioning trajectories to occur before a random short time, leading to refined estimates of hitting probabilities.} \\
{However, while previous work \cite{lawley2022extremehiting} has characterized extreme hitting probabilities in the limit of many searchers, our contribution provides a complementary boundary-layer analysis that explicitly captures the transition between competing asymptotic regimes. This allows us to derive closed-form approximations that depend directly on geometric parameters such as distances to targets and their relative configuration.} While analytical expressions for the distribution and scaling of the fastest arrival time are now well established, their direct use in simulations has remained limited. Most stochastic simulation methods (e.g., Gillespie algorithm or direct Brownian dynamics) become computationally prohibitive when extreme events must be resolved accurately, especially when large numbers of particles are involved.\\
{ The distribution of extreme first-passage times has been rigorously analyzed in \cite{lawley2020distribution}, who derived explicit formulas for the joint distribution of the first $k$ arrivals and their moments. These results provide an exact probabilistic characterization of extreme statistics and highlight the role of the Lambert $W$ function in describing limiting distributions. In contrast, the present work transform asymptotic first-passage results into an explicit computational recursive sampling algorithm that generates successive arrival times directly from asymptotic survival laws, avoiding the need to evaluate or invert closed-form distributions. This approach allows for efficient simulation in regimes where direct sampling from the exact distributions would be computationally impractical. These distributions are analytically involved, and their numerical inversion becomes challenging for large $n$, for joint distributions of successive arrivals, or when extensions such as multiple targets or time-dependent injection are considered. \\
Moreover, beyond the setting considered \cite{lawley2020distribution}, we incorporate multi-target geometries, time-dependent emission processes through multiple convolutions, and particle killing mechanisms within a single unified algorithm. These extensions  are essential for applications in spatially structured stochastic systems, where both geometry and temporal injection significantly influence extreme arrival statistics.} { Time-dependent emission profiles can be used to describe molecular release processes in biological systems, where particles are not introduced instantaneously but follow stochastic delayed production dynamics.} \\
{ We constructed here a recursive sampling framework based on conditional survival probabilities, allowing successive arrival times to be generated through local updates of the cumulative distribution. This approach avoids the explicit inversion of complex distributions and extends to incorporate splitting probabilities, time-dependent emission via convolution, and killing processes within a unified algorithmic structure. As a result, the method provides a trajectory-free and computationally efficient alternative to both direct Brownian simulations and inverse CDF sampling of known distributions.}\\
{ The present acceleration algorithm should be interpreted as a distribution-based approach, where the simulation is performed at the level of first-passage statistics rather than particle trajectories, in contrast to traditional Brownian or Monte Carlo methods. } Finally, our work here bridges this gap by providing an efficient simulation algorithm that bypasses the need for full trajectory generation. Building on asymptotic distributions for both arrival times and splitting probabilities, we demonstrate that it is possible to reproduce the extreme statistics of escape events through direct sampling, accelerating simulations by several orders of magnitude.\\
Moreover, the present work extends this framework in two critical directions: (1) we generalize to non-instantaneous injection profiles by incorporating time-dependent emission kernels into the extreme statistics formalism; and (2) we propose a recursive algorithm to simulate successive arrivals (first \(k\) particles), with and without exponential killing or time-dependent emission, thus capturing more complex event chains beyond the first escape.\\
{  while recent studies \cite{grebenkov2025, lawyley2025inhomogenous} have investigated the impact of time-dependent boundaries or trapping rates on survival probabilities, using limiting distributions and moments, they provided analytical insights into the role of temporal heterogeneity in extreme statistics. In contrast, we analyze here the interplay between structured emission and extreme arrival statistics by deriving explicit asymptotic formulas and integrating that framework into our algorithm. Indeed, rather than deriving only distributional properties, we constructed here  a recursive computational framework that enables the direct simulation of ordered arrival times under general emission profiles. In particular, we derive explicit asymptotic expressions for the mean fastest arrival time across multiple regimes and obtain inversion formulas involving the Lambert $W$ function, which allow efficient numerical evaluation.  To conclude, our approach integrates time-dependent injection with multi-target competition and particle killing within a single algorithmic structure, with applications to processes ubiquitous in cellular biology (e.g., calcium microdomains, signal transduction cascades, viral entry) but also in materials science, search theory, and communication networks, where first arrivals determine the system response.}
\subsection{code availability}
The code is now available at url \url{https://zenodo.org/records/19656299}
\section{Appendix}
\subsection{Appendix A : Short-time asymptotics for multiple targets}
For M small absorbing sites, each one is associated with a given survival probability:
\begin{equation*}
    S_i(t) = 1-f_i(t).
\end{equation*}
In the short-time approximation, there is independence between the absorbing targets and thus $f_i(t) \ll 1$. In first-order approximation, the total survival probability is
\begin{align*}
    S(t) = \prod_{i=1}^M S_i(t) = \prod_{i=1}^M(1-f_i(t)) \approx 1 -\sum_{i=1}^M f_i(t)
\end{align*}
\subsection{Appendix B : Conditional distribution of successive arrival time}
Using \ref{eq:order}, we have
\begin{align*}
    \mathbb{P}(\tau_{k+1}>t|\tau_k=t_k) &= \frac{\mathbb{P}(\tau_{k+1}>t, \tau_k = t_k)}{\mathbb{P}(\tau_k=t_k)} \\
    &= \frac{k \binom{n}{k} F(t_k)^{k-1} (1 - F(t))^{n - k} f(t_k)}{k \binom{n}{k} F(t_k)^{k-1} (1 - F(t_k))^{n - k} f(t_k)} \\
    &=\left(\frac{1-F(t)}{1-F(t_k)}\right)^{n-k} \\
    & = \exp((n-k)(\ln(1-F(t))-\ln(1-F(t_k)) \\
    &\approx \exp(-(n-k)(F(t)-F(t_k))
\end{align*}
\subsection{Appendix C : Asymptotic expressions for the splitting probabilities}
For dimension 1, we shall now compute  the integral $I(\lambda, n)$ asymptotically in the limit of large $n$.
\paragraph{For small $\lambda < 1$} we note:
\beq
    G(\lambda,n) = 1-I(\lambda,n) = \frac{1}{\sqrt{\pi}}\int_0^\infty \exp(-\frac{u}{\sqrt{\pi}})(1-\exp(-\frac{\lambda u^{\frac{1}{\lambda^2}}}{\sqrt{\pi}n^{\frac{1}{\lambda^2}-1}})) du
\eeq
As $\lambda < 1$, we have $\frac{1}{\lambda^2} - 1 > 0$, and therefore for sufficiently large $n$, provided that $u$ is bounded, we can expand the second exponential term around 0. Moreover, given the first exponential term, we can bound and neglect the remainder of the integral when $u$ is large. To do so, we split the integral at the point $A(n) = Cn^{1-\lambda^2 }$  with $C > 0$,
\begin{align*}
    G(\lambda,n) &= \frac{1}{\sqrt{\pi}}\int_0^{A(n)} \exp \left(-\frac{u}{\sqrt{\pi}} \right)\left(1-\exp \left(-\frac{\lambda u^{\frac{1}{\lambda^2}}}{\sqrt{\pi}n^{\frac{1}{\lambda^2}-1}}\right)\right) du \\
    & \ \ \ \ \ \ \ \ \ \ \  \ \ \ \ \ + \frac{1}{\sqrt{\pi}}\int_{A(n)}^\infty \exp\left(-\frac{u}{\sqrt{\pi}} \right)\left(1-\exp \left(-\frac{\lambda u^{\frac{1}{\lambda^2}}}{\sqrt{\pi}n^{\frac{1}{\lambda^2}-1}}\right)\right) du \\
    & = \frac{1}{\sqrt{\pi}}\int_0^{A(n)} \exp\left(-\frac{u}{\sqrt{\pi}} \right)\sum_{k=1}^\infty \frac{(-1)^{k+1}}{k!} \left(\frac{\lambda u^{\frac{1}{\lambda^2}}}{\sqrt{\pi}n^{\frac{1}{\lambda^2}-1}}\right)^k du + O\left(\int_{A(n)}^\infty \exp\left(-\frac{u}{\sqrt{\pi}} \right)  du\right) \\
    & =\sum_{k=1}^\infty \frac{(-1)^{k+1}}{k!} \left(\frac{\lambda }{\sqrt{\pi}n^{\frac{1}{\lambda^2}-1}}\right)^k \frac{1}{\sqrt{\pi}}\int_0^{A(n)}u^{\frac{k}{\lambda^2}}\exp \left( -\frac{u}{\sqrt{\pi}} \right)du + O\left( \exp \left(-\frac{Cn^{1-\lambda^2}}{\sqrt{\pi}} \right) \right) \\
    & = \sum_{k=1}^m \frac{(-1)^{k+1}}{k!} \left(\frac{\lambda }{\sqrt{\pi}n^{\frac{1}{\lambda^2}-1}}\right)^k \frac{1}{\sqrt{\pi}}\int_0^{\infty}u^{\frac{k}{\lambda^2}}\exp \left( -\frac{u}{\sqrt{\pi}} \right)du \\
    &\ \ \ \ \ \ \ \ \ \ \  \ \ \ \ \ -\sum_{k=0}^m \frac{(-1)^{k+1}}{k!} \left(\frac{\lambda }{\sqrt{\pi}n^{\frac{1}{\lambda^2}-1}}\right)^k \frac{1}{\sqrt{\pi}}\int_{A(n)}^{\infty}u^{\frac{k}{\lambda^2}}\exp \left( -\frac{u}{\sqrt{\pi}} \right)du \\
    &\ \ \ \ \ \ \ \ \ \ \  \ \ \ \ \ +\sum_{k=m+1}^\infty \frac{(-1)^{k+1}}{k!} \left(\frac{\lambda }{\sqrt{\pi}n^{\frac{1}{\lambda^2}-1}}\right)^k \frac{1}{\sqrt{\pi}}\int_0^{A(n)}u^{\frac{k}{\lambda^2}}\exp \left( -\frac{u}{\sqrt{\pi}} \right)du \\
    & \ \ \ \ \ \ \ \ \ \ \  \ \ \ \ \ + O\left( \exp \left(-\frac{n^{1-\lambda^2}}{\sqrt{\pi}} \right) \right) \ \ \text{avec} \ \ m =\left\lfloor \frac{C\lambda}{\sqrt{\pi}}\right\rfloor\\
\end{align*}
As the terms are exponentially small next to the leading term, we get :
\begin{align}
    I(\lambda,n) &\approx \sum_{k=0}^m \frac{(-1)^{k}}{k!} \Gamma\left(\frac{k}{\lambda^2}+1\right) \left( \frac{\sqrt{\pi}}{n} \right)^{k(\frac{1}{\lambda^2}-1)}  \lambda^k \\
    & = 1 - \lambda \Gamma\left(1+\frac{1}{\lambda^2}\right)\left(\frac{\sqrt{\pi}}{n} \right)^{\frac{1}{\lambda^2}-1}
\end{align}
\paragraph{For large $\lambda > 1$} the exponential term decays rapidly unless $u$ is small. The term in the second exponential tends to infinity much faster than the first one, so we restrict ourselves to considering the integration domain where $u$ is small and we split the integral at $u=C$ where $C$ is an arbitrary constant. We again expand the integrand:
\begin{align*}
I(\lambda,n) &= \frac{1}{\sqrt{\pi}} \int_0^C \sum_{k=0}^\infty \frac{(-1)^k u^k}{k! \sqrt{\pi}^k} \exp\left( - \frac{\lambda n^{1 - 1/\lambda^2} u^{1/\lambda^2}}{\sqrt{\pi}} \right) du + O\left( e^{- C^{1/\lambda^2} n^{1 - 1/\lambda^2}/\sqrt{\pi}} \right).
\end{align*}
Substituting $u = v / n^{\lambda^2 - 1}$ and expanding as before, we obtain:
\begin{align*}
    I(\lambda,n) &= \frac{1}{\sqrt{\pi}}\int_0^C \sum_{k=0}^\infty \frac{(-1)^k u^k}{k!\sqrt{\pi}^k} \exp\left( -\frac{\lambda n^{1-\frac{1}{\lambda^2}}u^{\frac{1}{\lambda^2}}}{\sqrt{\pi}} \right) du + O\left( \exp\left(-\frac{C^{\frac{1}{\lambda^2}}n^{1-\frac{1}{\lambda^2}}}{\sqrt{\pi}} \right) \right)\\
    &= \frac{1}{\sqrt{\pi}}\sum_{k=0}^\infty \int_0^C \frac{(-1)^k u^k}{k!\sqrt{\pi}^k} \exp\left( -\frac{\lambda n^{1-\frac{1}{\lambda^2}}u^{\frac{1}{\lambda^2}}}{\sqrt{\pi}} \right) du + O\left( \exp\left(-\frac{C^{\frac{1}{\lambda^2}}n^{1-\frac{1}{\lambda^2}}}{\sqrt{\pi}} \right) \right)\\
    &= \frac{1}{\sqrt{\pi}}\sum_{k=0}^\infty \int_0^{Cn^{\lambda^2-1}}\frac{(-1)^k u^k}{k!\sqrt{\pi}^k n^{(k+1)(\lambda^2-1)}}\exp\left(-\frac{\lambda u^{\frac{1}{\lambda^2}}}{\sqrt{\pi}} \right)du + O\left( \exp\left(-\frac{C^{\frac{1}{\lambda^2}}n^{1-\frac{1}{\lambda^2}}}{\sqrt{\pi}} \right) \right)\\
    &= \frac{1}{\sqrt{\pi}}\sum_{k=0}^m \int_0^{\infty}\frac{(-1)^k u^k}{k! \sqrt{\pi}^k n^{(k+1)(\lambda^2-1)}}\exp\left(-\frac{\lambda u^{\frac{1}{\lambda^2}}}{\sqrt{\pi}} \right)du \\
    &\qquad \qquad + \ \frac{1}{\sqrt{\pi}}\sum_{k=0}^m \int_{Cn^{\lambda^2-1}}^{\infty}\frac{(-1)^k u^k}{k!\sqrt{\pi}^k n^{(k+1)(\lambda^2-1)}}\exp\left(-\frac{\lambda u^{\frac{1}{\lambda^2}}}{\sqrt{\pi}} \right)du  \\
    &\qquad \qquad + \ \frac{1}{\sqrt{\pi}}\sum_{k=m+1}^\infty \int_{0}^{Cn^{\lambda^2-1}}\frac{(-1)^k u^k}{k!\sqrt{\pi}^k n^{(k+1)(\lambda^2-1)}}\exp\left(-\frac{\lambda u^{\frac{1}{\lambda^2}}}{\sqrt{\pi}} \right)du  \\
    & \qquad \qquad + \ O\left( \exp\left(-\frac{C^{\frac{1}{\lambda^2}}n^{1-\frac{1}{\lambda^2}}}{\sqrt{\pi}} \right) \right) \qquad \text{ with } \ m = \left\lfloor \frac{C}{\sqrt{\pi}} \right\rfloor.
\end{align*}
All those terms are exponentially small and thus we get :
\begin{align*}
    I(\lambda,n) &\approx \sum_{k=1}^{m+1} \frac{(-1)^{k+1}\Gamma(k\lambda^2+1)}{k!}\left( \frac{\sqrt{\pi}}{n}\right)^{k(\lambda^2-1)}\frac{1}{\lambda^{k\lambda^2}} \\
    &\approx \lambda\Gamma(\lambda^2)\left( \frac{\sqrt{\pi}}{n \lambda} \right)^{\lambda^2-1}
\end{align*}
\normalem
\bibliographystyle{ieeetr}
\bibliography{ref_general}
\end{document}